\documentclass[11pt,a4paper]{article}
\usepackage[utf8]{inputenc}
\usepackage{amscd,amsmath,amstext,amsfonts,amsbsy,amssymb,amsthm,eufrak}
\usepackage[left=2.5cm,right=2cm,top=2.5cm,bottom=2.5cm]{geometry}
\usepackage[unicode,colorlinks,linkcolor=blue,citecolor=red]{hyperref}
\usepackage{authblk}
\usepackage{xcolor}
\usepackage{graphicx}
\usepackage{multirow, dsfont}

\newtheorem{thm}{Theorem}[section]
\newtheorem{lem}[thm]{Lemma}
\newtheorem{proposition}[thm]{Proposition}

\newtheorem*{remark}{Remark}

\def\Z{\mathbb{Z}}
\def\N{\mathbb{N}}
\def\R{\mathbb{R}}
\def\L{\mathbb{L}}
\def\E{\mathbb{E}}
\def\P{\mathbb{P}}
\def\V{\mathcal{V}}
\def\eg{\textit{e.g.}\ }

\newcommand{\Var}{\operatorname{Var}}
\newcommand{\pen}{\operatorname{pen}}
\newcommand{\Crit}{\operatorname{Crit}}

\title{Adaptive Density Estimation Using Projection Kernels and Penalized Comparison to Overfitting}
\author{Van Ha Hoang \thanks{Faculty of Mathematics and Computer Science, Vietnam National University - Ho Chi Minh City, Viet Nam; \texttt{E-mail:\,hvha@hcmus.edu.vn}},
Tien Dat Nguyen \thanks{Faculty of Mathematics and Computer Science, Vietnam National University - Ho Chi Minh City, Viet Nam; \texttt{E-mail:\,ndat@hcmus.edu.vn}},
Thi Mong Ngoc Nguyen \thanks{Faculty of Mathematics and Computer Science, Vietnam National University - Ho Chi Minh City, Viet Nam; \texttt{E-mail:\,ngtmngoc@hcmus.edu.vn}}}

\begin{document}
\maketitle

\begin{abstract}
In this work, we study wavelet projection estimators for density estimation,
based on compactly supported $\mathcal S$-regular scaling functions. The main
issue is the choice of the resolution level, which determines the
bias--variance trade-off. We select this level by a Penalized Comparison to
Overfitting (PCO) criterion: each candidate estimator is compared in
$\L^2(\R)$ with a single overfitting reference, and the additional stochastic
fluctuation is corrected by an explicit penalty. For the selected estimator,
we prove a high-probability oracle inequality and a risk oracle inequality in
expectation. Over Besov balls of densities satisfying a common $\L^\infty$ bound, the risk
inequality yields the rate $n^{-2r/(2r+1)}$ uniformly over the class, with a
sufficient penalty threshold that can be chosen uniformly over the class. This rate has the
classical minimax order, and the procedure adapts to the unknown Besov
regularity. Numerical experiments on several density shapes confirm that the data-driven
level is close to the oracle one and adapts to the structure of the target.
\end{abstract}

\section{Introduction}
The selection of estimators has long been recognized as a central challenge in
nonparametric estimation, where the quality of an estimator depends on the
choice of certain tuning parameters. Such parameters typically take the form of
the bandwidth in kernel estimation, the resolution level in projection methods,
or the thresholding level in wavelet estimation. Determining these parameters in
a data-driven manner is essential for achieving adaptivity, that is, ensuring
that the estimator automatically adjusts to the unknown regularity of the
underlying function.

Several major approaches to adaptive estimation have been proposed in the
literature. One classical route involves the minimization of penalized empirical
risk, such as penalized likelihood or penalized least squares. Another widely
studied approach is based on pairwise comparison of estimators, with Lepski's
method \cite{lepskii1991problem} being a landmark contribution that provides an
adaptive solution for bandwidth selection in kernel estimation. In parallel,
Donoho and collaborators \cite{donoho1994ideal, donoho1995wavelet} initiated a
major line of research on wavelet thresholding, which has since become a
cornerstone for adaptive estimation in high-dimensional and irregular contexts.
Furthermore, the minimum complexity principle introduced by Barron and Cover
\cite{barron1991minimum}, later developed into the general framework of model
selection by Birg\'e and Massart \cite{barron1999risk, birge2007minimal}, has
provided deep insights into the connection between adaptivity and penalization.

More recently, attention has shifted towards designing fully data-driven
procedures for tuning parameter selection, with particular emphasis on bandwidth
and resolution choice. In kernel estimation, Lerasle \cite{lerasle2016optimal}
proposed an optimal penalty formulation for bandwidth selection via penalized
least squares. The Goldenshluger--Lepski (GL) method
\cite{goldenshluger2011bandwidth} further advanced this field by combining
pairwise comparison with an implicit bias-variance calibration. Lacour and
Massart \cite{lacour2016minimal} clarified the role of the so-called minimal
penalty within this framework, while Comte and Rebafka
\cite{comte2016nonparametric} pointed out some computational limitations of the
GL procedure. In response, Lacour, Massart, and Rivoirard
\cite{lacour2017estimator} introduced the {Penalized Comparison to
Overfitting} (PCO) method, which compares candidate estimators to a single
overfitting reference and corrects for the additional variance through an
explicit penalty. This approach has several appealing features: it is
computationally tractable, avoids the instability of cross-validation, and does
not rely on plug-in estimates of unknown quantities beyond a single tuning
constant. Theoretical guarantees such as oracle inequalities and minimax rates
have been established for PCO, and subsequent works
\cite{varet2020numerical, comte2020bandwidth} have extended its scope and
examined its numerical performance.

Despite these advances, the application of PCO to wavelet methods remains
relatively unexplored. Wavelets offer a powerful representation system for
functions and densities, combining localization in both time and frequency
domains. They are particularly well suited for capturing heterogeneous features
such as smooth trends, local irregularities, or boundary effects. Classical
histogram or kernel estimators are often less flexible in such scenarios, whereas
wavelet-based methods provide a natural multiresolution framework that adapts to
varying smoothness. The thresholding paradigm for wavelet coefficients has
already proved successful in signal processing and density estimation. However,
the specific problem of selecting the wavelet resolution level in a principled,
fully data-driven way via PCO has not been addressed in the literature. This
leaves a significant gap in adaptive estimation theory: while bandwidth selection
by PCO is now well understood, its extension to the discrete, multiresolution
setting is largely open.

In this work, we develop a PCO-based selection rule for wavelet projection
density estimators. Starting from a collection
$\{\widehat f_N:N\in\mathcal N\}$, we compare each candidate with the
overfitting estimator $\widehat f_{N_{\max}}$ and select the resolution level
by minimizing an explicit penalized criterion. For projection estimators in
$\L^2$, this PCO criterion turns out to be algebraically equivalent, up to an
additive constant independent of $N$, to a penalized least-squares criterion of
Birg\'e--Massart type \cite{barron1999risk,birge2007minimal}; the value of the
PCO viewpoint is not a different selector, but an explicit reading of the penalty
as a correction for the additional variance induced by the overfitting
comparison, together with a calibration of its constant.

Our theoretical contribution is twofold. First, we prove a high-probability oracle inequality, showing that the selected estimator is
close to the best estimator in the collection up to a multiplicative constant
and a negligible remainder. Second, under a suitable lower bound on the penalty
constant, we establish an oracle inequality in expectation with no remainder
depending on the maximal resolution. Combined with the usual Besov bias bound, this gives the rate $n^{-2r/(2r+1)}$ uniformly over the densities $f$ with
$\|f\|_{B^r_{2,q}}\le L$ and $\|f\|_\infty\le U$, for every $1\le q\le\infty$
and every $r$ below the regularity $\mathcal S$ of the scaling function. This rate has the classical minimax order for squared $\L^2$-risk over
Besov balls; see, \eg Gin\'e and Nickl \cite[Theorem~6.3.10]{book:Gine-RNickl.2016}. The sufficient lower bound on the penalty constant depends on the prescribed bound $U$, but not on $r$, $q$ or $L$, so the procedure does not require knowledge of these parameters. Numerical experiments
further demonstrate the estimator's performance across a variety of models with
different structural properties.

\medskip
\newpage 
The paper is organized as follows. In Section~\ref{sec:construction}, we
introduce the collection of wavelet projection estimators, state the conditions
on the scaling function, and record the Besov bias bound.
Section~\ref{sec:selection} presents the PCO-inspired selection rule for choosing
the resolution level. Theoretical guarantees, including the oracle inequality and
the rate of convergence, are established in Section~\ref{sec:results}.
Section~\ref{sec:simulation} reports numerical experiments that highlight the
practical performance of the method, and Section~\ref{sec:proofs} collects the
proofs.

\bigskip 
\noindent\textbf{Notation.} Unless stated
otherwise, $\|\cdot\|$ and $\langle\cdot,\cdot\rangle$ denote the norm and inner
product of $\L^2(\R)$, $P_N$ is the orthogonal projection onto $\V_N$,
$f_N:=P_N f$, and $a\lesssim b$ means $a\le C\,b$ for a constant $C$ independent
of $n$ and $N$.

\section{The collection of wavelet projection estimators}\label{sec:construction}

\subsection{Projection kernel and estimators}\label{sec:proj-estimator}

Let $\varphi$ be a compactly supported scaling function generating an
orthonormal multiresolution analysis $(\V_j)_{j\in\Z}$ of $\L^2(\R)$. For
$j,k\in\Z$ set
\[
\varphi_{j,k}(x):=2^{j/2}\varphi(2^{j}x-k),\qquad x\in\R,
\]
so that, by the defining properties of a multiresolution analysis, for every
$j\in\Z$ the family $\{\varphi_{j,k}:k\in\Z\}$ is an orthonormal basis of $\V_j$,
and the spaces are nested, $\V_j\subset \V_{j+1}$. For the statistical procedure we
restrict attention to nonnegative levels $j=N\in\N$, and we write
\[
\V_N=\overline{\operatorname{span}}\{\varphi_{N,k}:k\in\Z\},
\qquad N\in\N,
\]
with $\V_N\subset \V_{N+1}$. The orthogonal projection of $\L^2(\R)$ onto $\V_N$ is
the kernel operator with reproducing kernel
\[
K_N(x,y):=\sum_{k\in\Z}\varphi_{N,k}(x)\varphi_{N,k}(y),\qquad x,y\in\R,
\]
so that, for any $f\in\L^2(\R)$,
\begin{equation}\label{eq:proj-f}
f_N(x):=P_N f(x)=\int_\R K_N(x,y)f(y)\,dy
=\sum_{k\in\Z}\alpha_{N,k}\,\varphi_{N,k}(x),
\end{equation}
where the wavelet coefficients are
\[
\alpha_{N,k}:= \int_\R f(y)\varphi_{N,k}(y)dy = \E\bigl[\varphi_{N,k}(X_1)\bigr].
\]
Given the random sample $X_1, X_2, \ldots, X_n$, the natural unbiased estimator of $f_N$ replaces each coefficient
$\alpha_{N,k}$ by its empirical version:
\begin{equation}\label{eq:estimator}
\widehat f_N(x):=\sum_{k\in\Z}\widehat\alpha_{N,k}\,\varphi_{N,k}(x)
=\frac1n\sum_{u=1}^{n}K_N(X_u,x),
\end{equation}
where
\[
\widehat\alpha_{N,k}:=\frac1n\sum_{u=1}^{n}\varphi_{N,k}(X_u).
\]
Since $\E[\widehat\alpha_{N,k}]=\alpha_{N,k}$, we have
$\E[\widehat f_N]=f_N$. For a fixed sample and a fixed $N$, compact support of
$\varphi$ implies that only finitely many empirical coefficients
$\widehat\alpha_{N,k}$ are nonzero; moreover, for each evaluation point $x$, only
finitely many basis functions $\varphi_{N,k}(x)$ contribute. Thus
$\widehat f_N$ is exactly computable in practice.

Because $\V_N\subset \V_{N_{\max}}$, every candidate $\widehat f_N$ and the
overfitting reference $\widehat f_{N_{\max}}$ belong to the common space
$\V_{N_{\max}}$, and the relations needed for the PCO rule, in particular the
projection identity $\widehat f_N=P_N\widehat f_{N_{\max}}$ of
Section~\ref{sec:selection}, follow directly from orthogonality. 

\subsection{Technical conditions}\label{sec:tech-conditions}

The estimator~\eqref{eq:estimator} depends only on $\varphi$, and so do the
assumptions governing its behaviour. We impose the following two conditions on
the scaling function.

\begin{itemize}
\item[\textbf{(A1)}] $\varphi$ is compactly supported on $[-A,A]$ for some
integer $A\ge1$.
\item[\textbf{(A2)}] $\|\varphi\|_\infty\le\Phi_0$ for some constant
$\Phi_0>0$.
\end{itemize}

The next condition is used {only} to bound the bias over Besov classes
(Proposition~\ref{prop:bias}); it expresses that the projection kernel at level
$0$,
\[
K_0(x,y)=\sum_{k\in\Z}\varphi(x-k)\,\varphi(y-k),
\]
reproduces polynomials up to a sufficient order and has integrable off-diagonal
decay.

\begin{itemize}
\item[\textbf{(S)}] There exists an integer $\mathcal S\ge1$ such that
\begin{itemize}
\item[(i)] $\displaystyle\int_\R\varphi(u)\,du=1$, and for every $v\in\R$,
\[
\int_\R K_0(v,v+u)\,du=1,
\qquad
\int_\R K_0(v,v+u)\,u^{\ell}\,du=0,\quad \ell=1,\dots,\mathcal S-1;
\]
\item[(ii)] there exist constants $c_1,c_2\in(0,\infty)$ and a bounded
integrable function $\Phi:[0,\infty)\to\R_+$ such that
\[
\sup_{v\in\R}\bigl|K_0(v,v-u)\bigr|\le c_1\,\Phi(c_2|u|),
\qquad
\int_\R|u|^{\mathcal S}\,\Phi(|u|)\,du<\infty .
\]
\end{itemize}
\end{itemize}

Two remarks on these conditions are in order; they record the simplifications
relative to an earlier version of this work, in response to a referee.

\begin{remark}
Conditions \textbf{(A1)}--\textbf{(A2)} imply the property
\[
\sup_{x\in\R}\sum_{k\in\Z}\bigl|\varphi(x-k)\bigr|\le(2A+1)\,\Phi_0<\infty .
\]
Indeed, for each $x$ the support condition forces $\varphi(x-k)\ne0$ only for
the integers $k\in[x-A,x+A]$, of which there are at most $2A+1$, and each term
is bounded by $\Phi_0$. 
\end{remark}

\begin{remark}
The oracle inequalities (Theorem~\ref{thm:oracle} and the risk bound
Theorem~\ref{thm:rate}) use only \textbf{(A1)}--\textbf{(A2)}
together with $f\in\L^\infty$. Condition \textbf{(S)} enters exclusively
through the Besov bias bound (Proposition~\ref{prop:bias}), hence only in the
rate-of-convergence statement (Theorem~\ref{thm:rate}).
\end{remark}

Conditions \textbf{(A1)}--\textbf{(A2)} and \textbf{(S)} hold for the standard
compactly supported orthonormal wavelets. The Haar system covers low-regularity
cases under the appropriate Besov characterization, while sufficiently regular
compactly supported Daubechies wavelets cover higher smoothness regimes; the
exact admissible range of $r$ depends on the regularity and the number of
vanishing moments of the chosen wavelet system (the parameter $\mathcal S$
in~\textbf{(S)}, with $0<r<\mathcal S$). We refer to Gin\'e and Nickl
\cite[Section~4.3]{book:Gine-RNickl.2016} and H\"ardle \textit{et~al.}
\cite{book:Hardle-Kerkyacharian-Picard-Tsybakov:1998} for verification.
\medskip 

We close with the elementary pointwise kernel bound that controls 
variance terms in the sequel.

\begin{lem}\label{lem:kernel-bound}
Under \textbf{(A1)}--\textbf{(A2)}, for every $N\in\N$ and every $x\in\R$,
\begin{equation}\label{eq:kernel-bound}
\sum_{k\in\Z}\bigl|\varphi_{N,k}(x)\bigr|\le(2A+1)\,\Phi_0\,2^{N/2},
\qquad
\sum_{k\in\Z}\varphi_{N,k}^2(x)=K_N(x,x)\le(2A+1)\,\Phi_0^2\,2^{N}.
\end{equation}
In particular, at level $N=0$, $\sum_k|\varphi(x-k)|\le(2A+1)\Phi_0$ and
$\sum_k\varphi(x-k)^2\le(2A+1)\Phi_0^2$.
\end{lem}

\begin{proof}
Fix $x\in\R$ and $N\in\N$. The term $\varphi_{N,k}(x)=2^{N/2}\varphi(2^{N}x-k)$
is nonzero only when $|2^{N}x-k|\le A$, that is for
$k\in[2^{N}x-A,\,2^{N}x+A]\cap\Z$, a set of at most $2A+1$ integers; each such
term is bounded by $2^{N/2}\Phi_0$. Summing gives the first inequality
in~\eqref{eq:kernel-bound}. For the second, the same counting argument yields
$\sum_k\varphi_{N,k}^2(x)\le(2A+1)\cdot2^{N}\Phi_0^2$, and the identity
$\sum_k\varphi_{N,k}^2(x)=K_N(x,x)$ is the definition of the diagonal of the
projection kernel.
\end{proof}

\subsection{Besov classes and the bias bound}\label{sec:besov}

We measure the smoothness of the target density $f$ on the Besov scale
$B^{r}_{2,q}$, which controls the decay of the bias $\|f-f_N\|^2=\|f-P_N f\|^2$
as $N\to\infty$. We recall the wavelet characterization of these spaces. Since
Besov norms defined by sufficiently regular compactly supported wavelet bases
are equivalent, the basis $\{\varphi,\psi\}$ used to characterize $B^{r}_{2,q}$
may be taken to satisfy the regularity and approximation conditions below
without further restriction.

Let $\{\varphi,\psi\}$ be a compactly supported orthonormal wavelet basis of
regularity $\mathcal S>r$, with
$\varphi,\psi\in\mathcal C^{\mathcal S}(\R)$ and $D^{\mathcal S}\varphi,
D^{\mathcal S}\psi$ dominated by an integrable function. For $f\in\L^2(\R)$
write the detail coefficients
\[
\beta_{\ell,k}:=\int_\R f(y)\,\psi_{\ell,k}(y)\,dy,
\qquad
\psi_{\ell,k}(x):=2^{\ell/2}\psi(2^{\ell}x-k),
\quad\ell\in\N,\,k\in\Z,
\]
and set $\|\langle f,\psi_{\ell,\cdot}\rangle\|_2:=\bigl(\sum_{k}
|\beta_{\ell,k}|^2\bigr)^{1/2}$. For $0<r<\mathcal S$ and $1\le q\le\infty$ the
Besov space is
\[
B^{r}_{2,q}:=\bigl\{f\in\L^2(\R):\|f\|_{B^{r}_{2,q}}<\infty\bigr\},
\]
\[
\|f\|_{B^{r}_{2,q}}:=
\begin{cases}
\displaystyle
\|\langle f,\varphi_{\cdot}\rangle\|_2
+\Bigl(\sum_{\ell\ge0}2^{q\ell r}\,
\|\langle f,\psi_{\ell,\cdot}\rangle\|_2^{q}\Bigr)^{1/q},
& 1\le q<\infty,\\[10pt]
\displaystyle
\|\langle f,\varphi_{\cdot}\rangle\|_2
+\sup_{\ell\ge0}2^{\ell r}\,\|\langle f,\psi_{\ell,\cdot}\rangle\|_2,
& q=\infty .
\end{cases}
\]
We refer to \cite[Section~4.3]{book:Gine-RNickl.2016} for the details of these spaces.

\begin{proposition}[Gin\'e and Nickl, 2016]\label{prop:bias}
Suppose $\varphi$ satisfies \textbf{(A1)}, \textbf{(A2)} and \textbf{(S)} for
some integer $\mathcal S\ge1$. If $f\in B^{r}_{2,q}$ for some $0<r<\mathcal S$
and $1\le q\le\infty$, then there is a constant $C_{\mathrm{Besov}}<\infty$,
depending only on $\varphi$ (equivalently on $K_0$), such that for every $N\in\N$
\[
\bigl\|\E[\widehat f_N]-f\bigr\|^2=\|f_N-f\|^2
\le C_{\mathrm{Besov}}\,\|f\|_{B^{r}_{2,q}}^2\,2^{-2Nr}.
\]
\end{proposition}

This is the classical approximation bound for the projection $f_N=P_N f$. The
moment conditions in \textbf{(S)}(i) give the polynomial reproduction needed for
the order-$\mathcal S$ approximation, while the off-diagonal decay in
\textbf{(S)}(ii) controls the kernel away from the diagonal; together with the
wavelet characterization above they yield the stated decay. See, \eg
\cite[Proposition~4.3.8]{book:Gine-RNickl.2016} or
\cite[Section~9.5]{book:Hardle-Kerkyacharian-Picard-Tsybakov:1998}.

\section{Selection of the resolution level by the PCO method}\label{sec:selection}

The estimator $\widehat f_N$ of~\eqref{eq:estimator} depends on the resolution
level $N$, which governs the bias--variance trade-off and which we now choose in
a data-driven way. We select this level by adapting the Penalized Comparison to
Overfitting (PCO) method of Lacour, Massart and Rivoirard
\cite{lacour2017estimator}. Originally introduced for kernel density estimation,
PCO compares each candidate estimator with an overfitting benchmark and adds a
penalty to prevent overfitting. Here we transfer this principle to wavelet
projection estimators.

\medskip

Recall that $f$ denotes the unknown density and $\widehat f_N$ its wavelet
projection estimator at level $N$. Considering the classical bias--variance
decomposition
\begin{equation}\label{eq:bias-var}
\E\bigl[\|f-\widehat f_N\|^2\bigr]
=\bigl\|f-\E[\widehat f_N]\bigr\|^2
+\E\bigl[\|\E[\widehat f_N]-\widehat f_N\|^2\bigr]
=:B_N+V_N.
\end{equation}
This suggests considering a criterion of the form
\[
\Crit(N):=\widehat B_N+\widehat V_N,
\]
where $\widehat B_N$ and $\widehat V_N$ estimate the bias and variance terms,
respectively. Minimizing such a criterion is intended to mimic the minimization
of the $\L^2$-risk.

We begin with the variance term, which admits an exact expression. By
orthonormality of $\{\varphi_{N,k}\}_{k}$ and Parseval's identity,
\begin{align*}
V_N &= \E\big[ \| \widehat f_N - f_N \|^2\big]=\E\Bigl[\Bigl\|\sum_{k}(\widehat\alpha_{N,k}-\alpha_{N,k})
\varphi_{N,k}\Bigr\|^2\Bigr]
=\sum_{k}\Var(\widehat\alpha_{N,k}) \\
&=\frac1n\sum_{k}\Var\bigl(\varphi_{N,k}(X_1)\bigr)
\le\frac1n\,\E\Bigl[\sum_{k}\varphi_{N,k}^2(X_1)\Bigr].
\end{align*}

Lemma~\ref{lem:kernel-bound} bounds the inner sum pointwise by
$(2A+1)\Phi_0^2 2^N$, so that
\begin{equation}
	\label{eq:variance-bound}
	V_N
	\le
	(2A+1)\Phi_0^2\frac{2^N}{n}.
\end{equation}

{The upper bound~\eqref{eq:variance-bound} has the correct order in
$2^N/n$. Indeed, by orthonormality of $\{\varphi_{N,k}\}_k$,
\[
	V_N = \frac{1}{n}\sum_k \Var\!\bigl(\varphi_{N,k}(X_1)\bigr)
	= \frac{1}{n}\left(\E\bigl[K_N(X_1,X_1)\bigr] - \|f_N\|^2 \right).
\]
To see that the order $2^N/n$ cannot in general be improved, consider the density
\[
	f=\mathbf 1_{[0,1]}.
\]
Writing
\[
	K_N(x,x)= 2^N G(2^N x), \qquad G(u):=\sum_{k\in\mathbb Z}\varphi(u-k)^2,
\]
we note that $G$ is $1$-periodic and, by orthonormality,
\[
	\int_0^1 G(u)\,du=\int_{\mathbb R}\varphi(u)^2\,du=1.
\]
Hence
\[
	\E\bigl[K_N(X_1,X_1)\bigr]=\int_0^1 K_N(x,x)\,dx=2^N.
\]
Moreover, since $f_N=P_Nf$ and $P_N$ is an orthogonal projection,
\[
	\|f_N\|^2\le \|f\|^2=1.
\]
Therefore,
\[
	V_N=\frac{1}{n}\left(2^N-\|f_N\|^2	\right)\ge\frac{2^N-1}{n}.
\]
Together with~\eqref{eq:variance-bound}, this shows that, for this bounded density,
\[
	\frac{2^N-1}{n}\le V_N\le(2A+1)\Phi_0^2\frac{2^N}{n}.
\]
Thus the upper bound in~\eqref{eq:variance-bound} is sharp in order: the variance is naturally of order $2^N/n$ for projection estimators. Only the upper bound~\eqref{eq:variance-bound} is used in the construction of the criterion and in the oracle inequalities.
}

On the other hand, estimating the bias term is more difficult.
Recall that $N_{\max}$ is the highest resolution level under consideration. Since
$f_{N_{\max}}=\E\bigl(\widehat f_{N_{\max}}\bigr)$ provides a good approximation
of $f$, the squared difference $\|f_{N_{\max}}-f_N\|^2$ can be seen as a proxy
for the bias $B_N$.

A natural idea is to estimate this quantity empirically via
$\|\widehat f_{N_{\max}}-\widehat f_N\|^2$. However, this introduces an
additional variance term. Indeed, we have the decomposition
\begin{equation}\label{eq:hatdiff-decomp}
\widehat f_{N_{\max}}-\widehat f_N
=\bigl(\widehat f_{N_{\max}}-f_{N_{\max}}\bigr)
-\bigl(\widehat f_N-f_N\bigr)
+\bigl(f_{N_{\max}}-f_N\bigr),
\end{equation}
which yields the following identity:
\begin{equation}\label{eq:proxy-bias}
\E\bigl[\|\widehat f_{N_{\max}}-\widehat f_N\|^2\bigr]
=\|f_{N_{\max}}-f_N\|^2
+\E\bigl[\|\widehat f_{N_{\max}}-f_{N_{\max}}-(\widehat f_N-f_N)\|^2\bigr].
\end{equation}

We will control the fluctuation term on the right of~\eqref{eq:proxy-bias}. The key is the following projection
identity, which expresses that $\widehat f_N$ is itself the projection of the
overfitting reference onto $\V_N$.

\begin{lem}\label{lem:proj-identity}
For every $0\le N\le N_{\max}$,
\begin{equation}\label{eq:proj-identity}
\widehat f_N=P_N\widehat f_{N_{\max}},
\qquad\text{and likewise}\qquad
f_N=P_N f_{N_{\max}} .
\end{equation}
In particular, with $I$ the identity operator on $\L^2(\R)$,
$\widehat f_{N_{\max}}-\widehat f_N=(I-P_N)\widehat f_{N_{\max}}$ and
$f_{N_{\max}}-f_N=(I-P_N)f_{N_{\max}}$.
\end{lem}

\begin{proof}
Both $\{\varphi_{N,j}\}_j$ and $\{\varphi_{N_{\max},k}\}_k$ are orthonormal
bases of $\V_N\subset \V_{N_{\max}}$ and $\V_{N_{\max}}$, so each $\varphi_{N,j}$
is reproduced by $P_{N_{\max}}$. The $j$-th coefficient of
$P_N\widehat f_{N_{\max}}$ in the basis $\{\varphi_{N,j}\}_j$ is
\begin{align*}
\langle P_N\widehat f_{N_{\max}},\varphi_{N,j}\rangle
& =\langle\widehat f_{N_{\max}},\varphi_{N,j}\rangle
=\frac1n\sum_{u}\sum_{k}\varphi_{N_{\max},k}(X_u)
\langle\varphi_{N_{\max},k},\varphi_{N,j}\rangle \\
& =\frac1n\sum_{u}\varphi_{N,j}(X_u)
=\widehat\alpha_{N,j},
\end{align*}
the third equality using $P_{N_{\max}}\varphi_{N,j}=\varphi_{N,j}$ (the sums are
finite by \textbf{(A1)}). This proves $\widehat f_N=P_N\widehat f_{N_{\max}}$.
For $f_N$, since $f_{N_{\max}}=P_{N_{\max}}f$ and $P_N P_{N_{\max}}=P_N$, we get
$P_N f_{N_{\max}}=P_N f=f_N$. The two displayed consequences are immediate.
\end{proof}

Set $D:=\widehat f_{N_{\max}}-f_{N_{\max}}\in \V_{N_{\max}}$. By applying Lemma~\ref{lem:proj-identity}, 
\[
\widehat f_{N_{\max}}-f_{N_{\max}}-(\widehat f_N-f_N)
=(I-P_N)\bigl(\widehat f_{N_{\max}}-f_{N_{\max}}\bigr)=(I-P_N)D .
\]

Since $P_N$ is an orthogonal projection, Pythagoras gives
\[
\|(I-P_N)D\|^2=\|D\|^2-\|P_N D\|^2 = \|\widehat f_{N_{\max}}-f_{N_{\max}}\|^2 - \| \widehat f_N-f_N\|^2 
\]
since $P_N D= P_N(\widehat f_{N_{\max}}-f_{N_{\max}}) = P_N \widehat f_{N_{\max}} - P_N  f_{N_{\max}}  = \widehat f_N - f_N$ by Lemma~\ref{lem:proj-identity}. \medskip 

Taking expectations yields
\begin{align}
\E\bigl[\|\widehat f_{N_{\max}}-f_{N_{\max}}-(\widehat f_N-f_N)\|^2\bigr] & = \E\big[ \|\widehat f_{N_{\max}}-f_{N_{\max}}\|^2\big] - \E\big[ \| \widehat f_N-f_N\|^2\big] \nonumber \\
& = V_{{N_{\max}}} - V_N. \nonumber 
\end{align}
{
Indeed, using the exact expression for $V_N$ above and the nestedness
of the projection spaces,
\[
V_{N_{\max}}-V_N
\le \frac{1}{n}\mathbb{E}\!\left[
K_{N_{\max}}(X_1,X_1)-K_N(X_1,X_1)\right].
\]
Moreover, for an index $j$, since $P_{j+1}-P_j$ is an orthogonal projection,
\[
0\le K_{j+1}(x,x)-K_j(x,x)
\le K_{j+1}(x,x)
\le (2A+1)\Phi_0^2 2^{j+1}.
\]
Hence, for any $N<N_{\max}$,
\begin{equation}\label{eq:bound-fluctuation}
	V_{N_{\max}}-V_N
	\le
	\frac{(2A+1)\Phi_0^2}{n}
	\sum_{j=N}^{N_{\max}-1}2^{j+1}
	=
	\frac{2(2A+1)\Phi_0^2}{n}
	\left(2^{N_{\max}}-2^N\right).
\end{equation}
When $N=N_{\max}$, the difference is trivially equal to zero.
Hence, the bound above is valid for every $N\le N_{\max}$.\\
The derivations above lead to the following fully data-driven model
selection criterion, which balances the approximation error and the
stochastic fluctuations. In this formulation, the additional numerical
factor $2$ only affects the numerical calibration of the penalty and
does not change its functional form in $N$, while terms depending solely
on $N_{\max}$ are irrelevant to the minimization over $N$.
}

For each resolution level $N\in\mathcal N$, we define the criterion

\begin{equation}\label{eq:Crit}
\Crit(N):=\bigl\|\widehat f_{N_{\max}}-\widehat f_N\bigr\|^2+\pen_\lambda(N),
\end{equation}
with penalty
\begin{equation}\label{eq:pen}
\pen_\lambda(N)
:=\frac{(2A+1)\,\Phi_0^2}{n}
\Bigl[\lambda\,2^{N}-\bigl(2^{N_{\max}}-2^{N}\bigr)\Bigr]
=\frac{(2A+1)\,\Phi_0^2}{n}
\Bigl[(\lambda+1)\,2^{N}-2^{N_{\max}}\Bigr].
\end{equation}

This criterion is motivated by the fact that $\bigl\|\widehat f_{N_{\max}}-\widehat f_N\bigr\|^2$ serves as a noisy estimate of the squared bias  $\bigl\| f_{N_{\max}}- f_N\bigr\|^2$ , but also includes a variance term that must be corrected. The penalty compensates for this overestimation through a correction
of the same order as the variance difference, together with a calibrated
penalty proportional to the variance of $\widehat f_N$.


\medskip 
The data-driven estimator $\widehat f_{\widehat N}$ is then defined by
\begin{equation}\label{eq:hatN}
\widehat N:=\operatorname*{arg\,min}_{N\in\mathcal N}\Crit(N)
=\operatorname*{arg\,min}_{N\in\mathcal N}
\Bigl\{\bigl\|\widehat f_{N_{\max}}-\widehat f_N\bigr\|^2+\pen_\lambda(N)\Bigr\},
\end{equation}
where the index set $\mathcal N$ will be detailed later.

\begin{remark}
For projection estimators in $\L^2$, the criterion~\eqref{eq:Crit} reduces to a
penalized least-squares criterion of Birg\'e--Massart type. Indeed, expanding the
square and using $\langle\widehat f_{N_{\max}},\widehat f_N\rangle
=\|\widehat f_N\|^2$ (which follows from the projection
identity~\eqref{eq:proj-identity}, since $\widehat f_N=P_N\widehat f_{N_{\max}}$),
\[
\bigl\|\widehat f_{N_{\max}}-\widehat f_N\bigr\|^2
=\|\widehat f_{N_{\max}}\|^2-\|\widehat f_N\|^2 .
\]
The term $\|\widehat f_{N_{\max}}\|^2$ does not depend on $N$. Moreover, for the
projection estimator one has the exact identity
$\tfrac1n\sum_{u=1}^{n}\widehat f_N(X_u)=\|\widehat f_N\|^2$, obtained from
$\widehat f_N=\sum_k\widehat\alpha_{N,k}\varphi_{N,k}$ and
$\widehat\alpha_{N,k}=\tfrac1n\sum_u\varphi_{N,k}(X_u)$ by orthonormality. Hence,
with the least-squares contrast
$\gamma_n(t):=\|t\|^2-\tfrac2n\sum_{u=1}^{n}t(X_u)$, we get
$\gamma_n(\widehat f_N)=-\|\widehat f_N\|^2$, and minimizing~\eqref{eq:Crit} over
$N\in\mathcal N$ is equivalent to minimizing
\[
\gamma_n(\widehat f_N)+\frac{(\lambda+1)(2A+1)\Phi_0^2}{n}\,2^{N},
\]
that is, a penalized least-squares criterion whose penalty is proportional to the
effective dimension $2^{N}$, exactly as in the model selection theory of Birg\'e
and Massart \cite{barron1999risk,birge2007minimal}. The value of the PCO
derivation is therefore not a new selector, but an explicit reading of the
penalty as a correction for the additional variance produced by the overfitting
comparison, together with a calibration of its constant.
\end{remark}

\newpage 
It remains to fix the collection $\mathcal N$ of candidate levels and the
overfitting level $N_{\max}$. The variance bound~\eqref{eq:variance-bound}
suggests that the stochastic size of $\widehat f_N$ is governed by an effective
complexity of order $2^{N}$; we therefore restrict the maximal level so that
this effective complexity remains below the sample size, up to logarithmic
factors. The relevant threshold, familiar from
\cite{lacour2017estimator,birge2007minimal}, is of order $n/(\ln n)^2$. We
therefore take
\begin{equation}\label{eq:N-set}
\mathcal N:=\Bigl\{N\in\N:\ N\le\bigl\lfloor\log_2\bigl(n/(\ln n)^2\bigr)
\bigr\rfloor\Bigr\},
\qquad
N_{\max}:=\max\mathcal N=\bigl\lfloor\log_2\bigl(n/(\ln n)^2\bigr)\bigr\rfloor,
\end{equation}
so that $2^{N_{\max}}\asymp n/(\ln n)^2$.

\section{Main results}\label{sec:results}
We now state the two theoretical guarantees satisfied by the data-driven
estimator $\widehat f_{\widehat N}$ of~\eqref{eq:hatN}: a
high-probability oracle-type comparison inequality, valid for any bounded
density, and the resulting rate of convergence
over Besov classes. Both are proved in Section~\ref{sec:proofs}. In what follows, $\mathcal N$ and $N_{\max}$ are
as in~\eqref{eq:N-set}, $c_\star(p)=c_1(p+1)$ is the constant of
Lemma~\ref{lem:good-event}, $c_{3}(p)$ is the constant of Lemma~\ref{lem:mixedcross}, and we recall
$V_2(n,N)=(2A+1)\Phi_0^2\,2^{N}/n+\ln n/n$ from~\eqref{eq:V2-def}.
The constant $c_\star(p)$ is explicit. Indeed, by	Proposition~\ref{prop:concentration}, one may take $c_1(p)=\frac{200}{9}(p+1)\bigl(1+\|f\|_\infty\bigr)$. Hence, by~\eqref{eq:cstar},
$c_\star(p) = c_1(p+1) = \frac{200}{9}(p+2)\bigl(1+\|f\|_\infty\bigr)$. Thus one may take the universal numerical constant $\kappa=200/9$. In particular, the dependence of $c_\star(p)$ on $\|f\|_\infty$ is explicit and linear. Consequently, the penalty threshold $\lambda_0(p,U)$ of Theorem~\ref{thm:rate}, defined just before that theorem, is of order $(p+2)(1+U)$, where $U$ is a prescribed upper bound on $\|f\|_\infty$. The second constant $c_3(p)$, coming from Lemma~\ref{lem:mixedcross}, enters only the negligible $\ln n/n$ remainder and plays no role in the leading trade-off.

\begin{thm}\label{thm:oracle}
{ 
Assume \emph{(A1)--(A2)} and $f\in\L^\infty(\R)$. Let $\theta\in(0,1)$, $p\ge1$ and $\lambda > (4+2\theta ) c_\star(p) - 1 $. Then, for $n$ sufficiently large, with probability at least
$1-2n^{-p}$, 
\[
\bigl\|\widehat f_{\widehat N}-f\bigr\|^2
\le\frac{C_{\theta,p, \lambda}}{1-\theta}\,
\min_{N\in\mathcal N} \Big( \bigl\|\widehat f_N-f\bigr\|^2  + \dfrac{2^{N}}{n} \Big)  
+
\frac{ 2 c_{\star}(p) + 4 \frac{c_{3}(p)}{\theta} }{1-\theta} \frac{\ln n}{n},   
\]
where $C_{\theta,p, \lambda} := \max \big\{ 1 + \theta ;  \bigl[(\lambda+1)+ 2\theta c_\star(p)\bigr]
(2A+1) \Phi_{0}^{2} \big\}$. 
}   
\end{thm}
{
Theorem~\ref{thm:oracle} is a high-probability oracle inequality.
On an event of probability at least $1-2n^{-p}$, the squared $\L^2$
error of the selected estimator is controlled by the best
error--complexity trade-off available in the collection, up to a
multiplicative constant and an additive remainder. The upper bound has
two parts.
\begin{itemize}
\item The leading term $	\frac{C_{\theta,p,\lambda}}{1-\theta} \min_{N\in\mathcal N}\left(\|\widehat f_N-f\|^2+\frac{2^N}{n}\right)$ shows that $\widehat f_{\widehat N}$ performs, up to a constant, as well as the estimator in the collection achieving the best trade-off between the estimation error $\|\widehat f_N-f\|^2$ and the variance-complexity
term $2^N/n$.
\item The remainder is of order $\ln n/n$, does not involve $N_{\max}$, and is negligible at the minimax oracle scale considered below.
\end{itemize}
The price for this sharper remainder is the condition $\lambda>(4+2\theta)c_\star(p)-1$
on the penalty constant. This condition allows the variance contribution
at the selected level $\widehat N$ to be absorbed by the penalty at that
level, rather than appearing in the remainder through a term of order
$2^{N_{\max}}/n$. In this respect, Theorem~\ref{thm:oracle} already
contains the key mechanism used in the risk bound below.
Theorem~\ref{thm:rate} then takes expectations and combines this oracle
trade-off with the Besov bias bound to obtain a rate of convergence
that is uniform over bounded Besov classes of densities.
}

\medskip 
We now turn to a risk bound in expectation and to the resulting
rate of convergence over bounded Besov classes of densities. Let $U>0$ be a
prescribed upper bound on $\|f\|_\infty$. For $0<r<\mathcal S$,
$1\le q\le\infty$, and $L>0$, we consider the class
\[
\mathcal D_{r,q}(L,U) := \Bigl\{ f:\R\to[0,\infty):
\int_\R f(x)\,dx=1,\quad \|f\|_{B^r_{2,q}}\le L,\quad \|f\|_\infty\le U \Bigr\}.
\]
For fixed $p$, the sufficient penalty threshold depends on the class
parameters only through $U$. We set
\[
c_\star(p,U) := \frac{200}{9}(p+2)(1+U), \qquad \lambda_0(p,U) := 5c_\star(p,U)-1.
\]
Since $c_\star(p)=\frac{200}{9}(p+2)(1+\|f\|_\infty)$ (see above), we have
$c_\star(p)\le c_\star(p,U)$ whenever $\|f\|_\infty\le U$. Hence every
$\lambda\ge\lambda_0(p,U)$ is admissible simultaneously for all densities satisfying
$\|f\|_\infty\le U$.

\begin{thm}[Risk bound and rate of convergence]\label{thm:rate}
Assume \emph{(A1)--(A2)--(S)}. Let $p\ge2$, $U>0$ and
$\lambda\ge\lambda_0(p,U)$. There exist an integer $n_0$, depending only on
$p$, and constants $C_1$, depending only on $\lambda$, and $C_2$, depending only
on $p$, $U$ and $\varphi$, such that for every probability density $f$ with
$\|f\|_\infty\le U$, every $n\ge n_0$ and every $N\in\mathcal N$,
\begin{equation}\label{eq:risk-oracle}
	\E\bigl\|\widehat f_{\widehat N}-f\bigr\|^2
	\;\le\; C_1\Bigl(\|f_N-f\|^2+(2A+1)\Phi_0^2\,\frac{2^{N}}{n}\Bigr)
	\;+\;C_2\,\frac{\ln n}{n}.
\end{equation}
In particular, $C_1$ and $C_2$ do not depend on $f$, $n$, $N$ or $N_{\max}$.
Consequently, for every $0<r<\mathcal S$, $1\le q\le\infty$ and $L>0$, there
exist an integer $n_1\ge n_0$, depending only on $p$ and $r$, and a constant
$C$, depending only on $\lambda$, $p$, $U$, $L$, $r$ and $\varphi$, such that
for all $n\ge n_1$,
\begin{equation}\label{eq:rate}
	\sup_{f\in\mathcal D_{r,q}(L,U)}\E\bigl\|\widehat f_{\widehat N}-f\bigr\|^2
	\;\le\; C\,n^{-\frac{2r}{2r+1}}.
\end{equation}
\end{thm}

\medskip
The risk oracle inequality~\eqref{eq:risk-oracle} is a bias--variance oracle
bound. Its leading term is the deterministic trade-off
\[
\min_{N\in\mathcal N}
\left\{
\|f_N-f\|^2+(2A+1)\Phi_0^2\frac{2^N}{n}
\right\},
\]
which is balanced at the oracle level $N_*$ satisfying
$2^{N_*}\asymp n^{1/(2r+1)}$, while the remainder term $C_2\ln n/n$ is of
smaller order. Since the threshold $\lambda_0(p,U)$ does not depend on $r$, $q$ or $L$, the
same penalty choice $\lambda\ge\lambda_0(p,U)$ is valid for every
$0<r<\mathcal S$, $1\le q\le\infty$ and $L>0$. Thus, once $U$ is fixed, the
procedure requires no knowledge of the smoothness parameters $r$ and $q$, nor
of the Besov radius $L$. The rate $n^{-2r/(2r+1)}$ in~\eqref{eq:rate} is the
classical minimax rate for density estimation over Besov balls; see, \eg
\cite[Theorem~6.3.10]{book:Gine-RNickl.2016}.

\section{Numerical experiments}\label{sec:simulation}

To illustrate the performance of the proposed adaptive estimator \(\widehat f_{\widehat N}\), built from~\eqref{eq:estimator} and selected by~\eqref{eq:hatN}, we conduct a Monte Carlo simulation study across several density models.

\subsection{Simulation setting}

We generate an i.i.d. sample \( X_1, \dots, X_n \) from four distinct probability distributions, chosen to represent different structural and smoothness characteristics:

\begin{itemize}
\item[\textbf{(M1)}] Standard Gaussian $X \sim \mathcal{N}(0,1)$.  This serves as a standard benchmark for symmetric, unimodal, and smooth densities.
\item[\textbf{(M2)}] Gaussian mixture $X\sim 0.25\,\mathcal{N}(0,1)+0.75\,\mathcal{N}(10,4)$. This model introduces multimodality and heteroscedasticity.

\item[\textbf{(M3)}] $X\sim \mathrm{Beta}(2,5)$: a bounded, asymmetric density on $[0,1]$ with a boundary effect at $0$.

\item[\textbf{{(M4)}}] {Claw density (see \textit{e.g.} Marron and Wand~\cite{marron1992exact}):
$X\sim \tfrac12\,\mathcal N(0,1)+\sum_{l=0}^{4}\tfrac{1}{10}\,\mathcal N\!\big(l/2-1,\,(1/10)^2\big)$. This benchmark consists of a broad Gaussian background with five sharp, narrow Gaussian spikes and is included to assess adaptation to fine local
structure.}
\end{itemize}

\begin{figure}[htbp]
\centering
\includegraphics[width=0.9\textwidth]{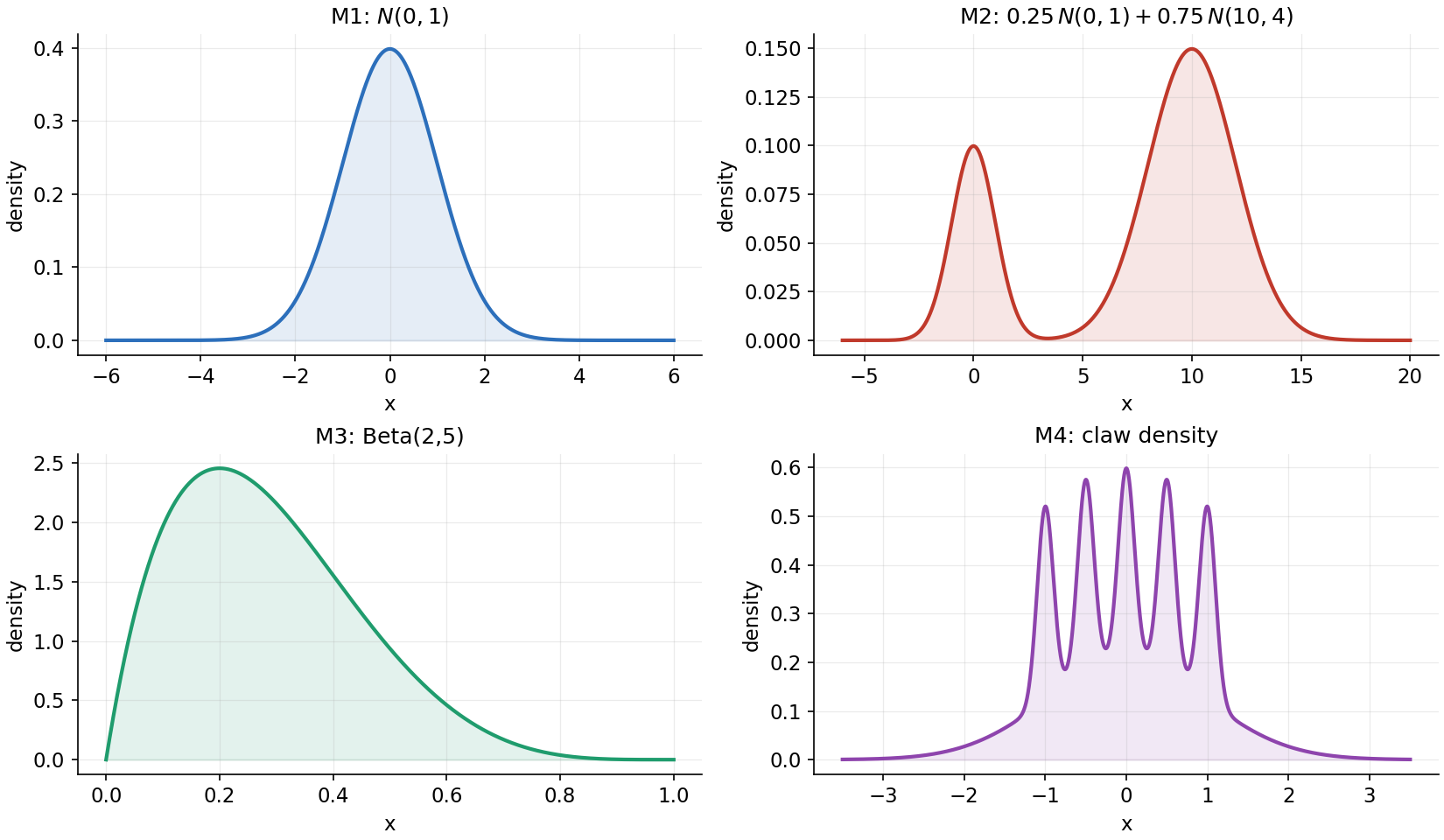}
\caption{True densities of the {four} simulation models {: M1 and M2 (top row), M3 and M4 (bottom row)}.}
\label{fig:true-densities}
\end{figure}

The true densities of the {four} models are shown in Figure~\ref{fig:true-densities}.

We use \(M=100\) Monte Carlo replications for all four models. The sample size is \(n=2^{12}\) for M1–M3 and \(n=2^{15}\) for M4. Errors are reported as
the Monte Carlo estimate of the mean integrated squared error (MISE).
\[
\mathrm{MISE}(\widehat f_N)=\E\|\widehat f_N-f\|^2\approx
\frac1M\sum_{j=1}^M \big\|\widehat f_N^{(j)}-f\big\|^2,
\]
where $\widehat f_N^{(j)}$ denotes the estimator computed from the $j$-th Monte Carlo sample, $j=1,2,\ldots,M$. 

We use the Daubechies scaling function with eight vanishing moments (\texttt{db8}). 
We have $\mathrm{supp}\,\varphi=[0,15]$
for \texttt{db8}, giving the half-width $A=15$, $\Phi_0=\|\varphi\|_\infty\approx
1.010$, and hence the penalty constant
\[
C_\varphi=(2A+1)\Phi_0^2\approx 31.65 .
\]
Numerically $\int\varphi=1$ and $\int\varphi^2=1$ to machine precision, as
required. The candidate set is determined by the sample size through
$N_{\max}=\big\lfloor\log_2\!\big(n/(\ln n)^2\big)\big\rfloor$, so for $n=2^{12}$
one obtains $N_{\max}=5$ and $\mathcal N=\{1,2,3,4,5\}$.
{The finely structured claw density M4 requires a finer resolution
than the other three; to give the candidate set enough room we use a larger
sample size $n=2^{15}$ for M4, which yields $N_{\max}=8$.}

\medskip 
We compare $\widehat f_{\widehat N}$ denoted by {WDE-PCO}
with three references:
\begin{itemize}
\item[(i)]  the {WDE-Oracle} $\widehat f_{N_{\mathrm{Oracle}}}$ with
$N_{\mathrm{Oracle}}=\arg\min_{N\in\mathcal N}\|\widehat f_N-f\|^2$, which uses
the unknown $f$ and is the best achievable level in the family;

\item[(ii)] a Gaussian kernel estimator with least-squares cross-validated bandwidth
({KDE-CV}); and

\item[(iii)]  the same with Silverman's robust rule-of-thumb bandwidth
({KDE-RoT}). 
\end{itemize}

\subsection{Calibration of $\lambda$}
\begin{figure}[htbp]
\centering
\includegraphics[width = 0.49\textwidth]{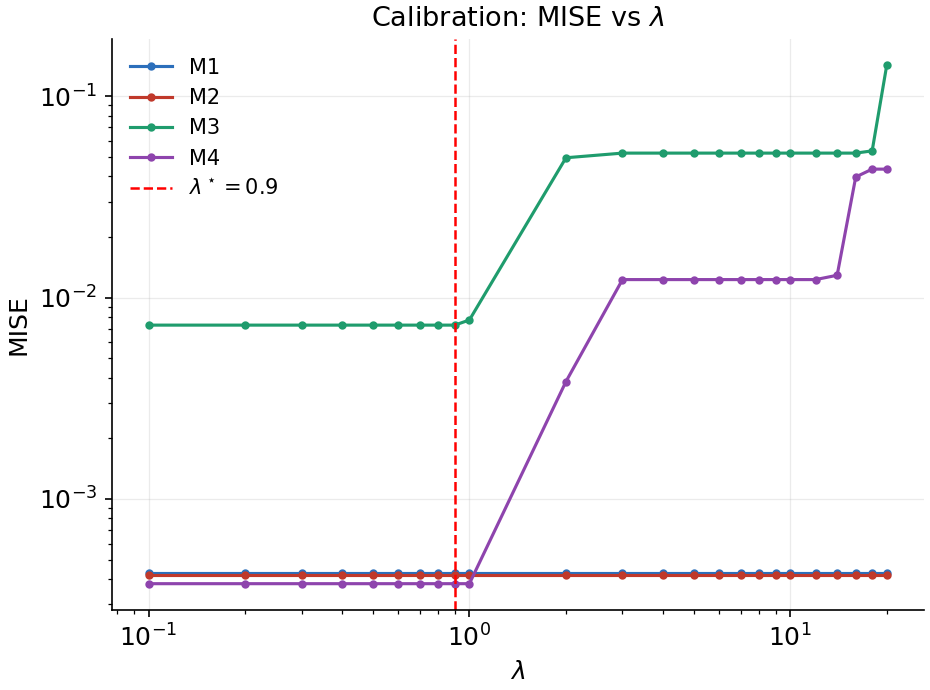}
\includegraphics[width=0.49\textwidth]{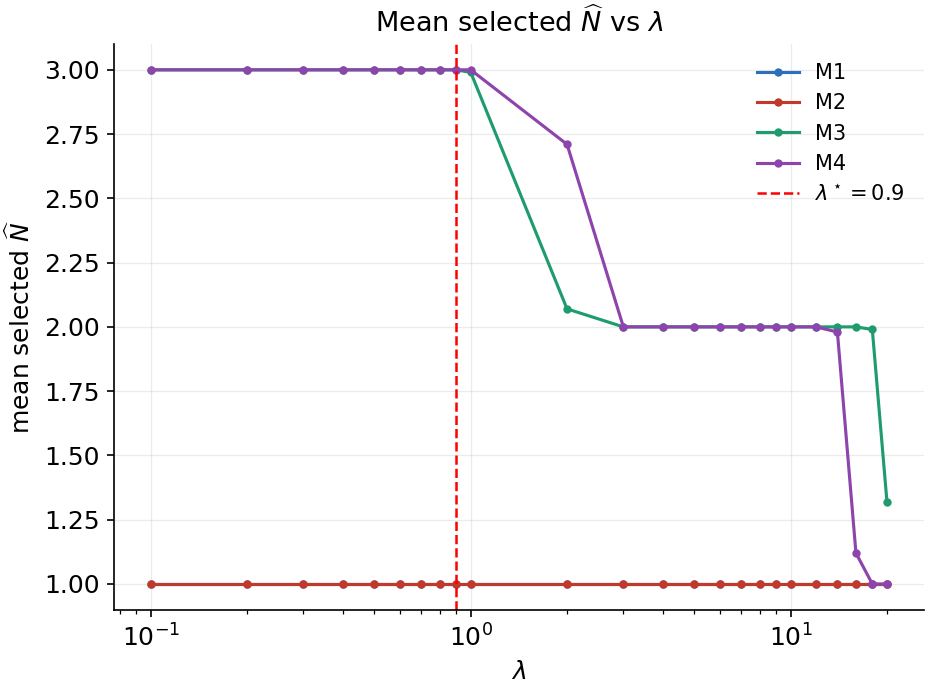}
\caption{Calibration of $\lambda$ for {all four models}. Left: $\mathrm{MISE}(\lambda)$ per model (log
scale; M1 and M2 are both $\approx 4\times10^{-4}$, as is M4 at its larger sample size); the dashed line marks
$\lambda^\star=0.9$. Right: mean selected level $\widehat N$.}
\label{fig:calibration}
\end{figure}

For each model and each $\lambda$ in the grid $\Lambda=\{0.1,0.2,\ldots,1,2,\ldots,10,12,14,16,18,20\}$,
we computed $\operatorname{MISE}_i(\lambda)$ over the $M=100$
replications. Figure~\ref{fig:calibration} shows a broad low-error region
for small values of $\lambda$. In particular, $\lambda^\star=0.9$ lies in
this stable region for all four models, and we therefore use this common
value in the comparisons below.

The calibration curves behave differently across the models. For M1 and M2, the chosen level stays at $\widehat N=1$ for every value of $\lambda$. These two densities are smooth and spread out, so the coarsest level already gives the best fit, and changing $\lambda$ makes no difference. For M3 the choice of $\lambda$ matters: the method selects $\widehat N=3$ when $\lambda\le 0.9$, then drops to $\widehat N=2$ once $\lambda \ge 2$, which makes the error much larger. {The finely structured claw density M4 behaves similarly to M3: it selects $\widehat N=3$ for $\lambda\le0.9$, and the error grows
as $\lambda$ increases further.} In all cases, the method never selects the largest level $N_{\max}$. 

{
The value $\lambda^\star=0.9$ selected in the calibration experiment is
smaller than the sufficient threshold $\lambda_0(p,U)=5c_\star(p,U)-1$ required by Theorem~\ref{thm:rate}, where $U$ is a prescribed upper bound on $\|f\|_\infty$. This does not contradict the theoretical result: $\lambda_0(p,U)$ is a
sufficient threshold obtained from worst-case concentration bounds and is
not intended to be necessary or numerically sharp. Such conservativeness
is common in penalized model-selection procedures and has also been
observed for PCO-type methods; see, e.g., \cite{lacour2017estimator,varet2020numerical}.\\
Although the theory does not cover the regime $\lambda<\lambda_0(p,U)$,
Figure~\ref{fig:calibration} shows empirically that the selection rule
remains stable around $\lambda^\star=0.9$. In particular, the maximal
resolution $N_{\max}$ is never selected. The smooth models M1 and M2
remain at the coarsest level $N=1$, whereas the more structured models
M3 and M4 select the higher level $\widehat N=3$ for small values of
$\lambda$ and move toward lower resolutions as the penalty increases.
}

\subsection{Simulation results}

\begin{figure}[htbp]
\centering
\includegraphics[width=\textwidth]{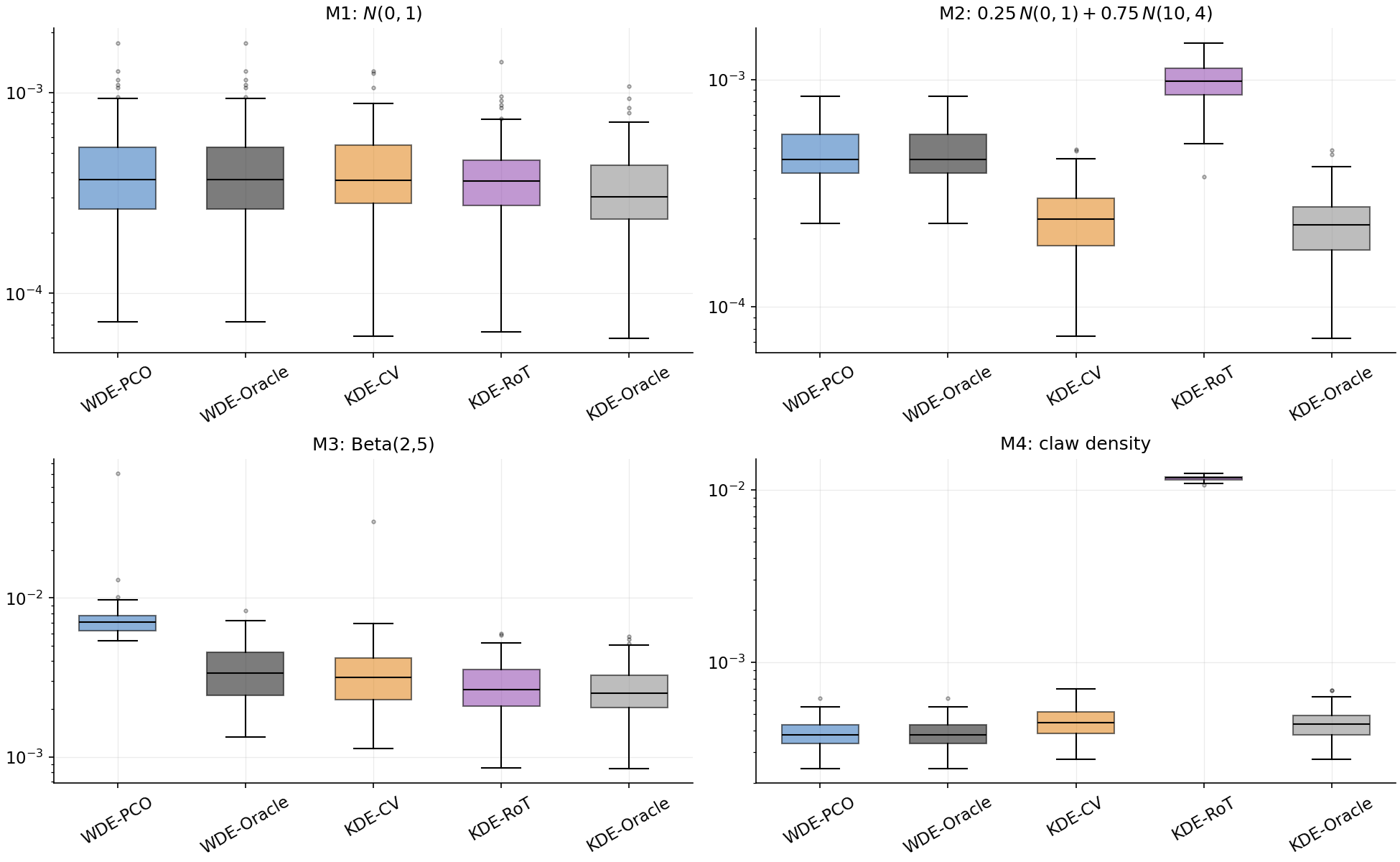}
\caption{Errors by method for each model: M1 and M2 (top row), M3
and M4 (bottom row). WDE-PCO coincides with WDE-Oracle on M1 and M2 and is nearly
indistinguishable from it on M4; note the shifted vertical scale on M4,
where the risks are those of the larger sample size $n=2^{15}$.}
\label{fig:boxplots}
\end{figure}

\begin{figure}[htbp]
\centering
\includegraphics[width=\textwidth]{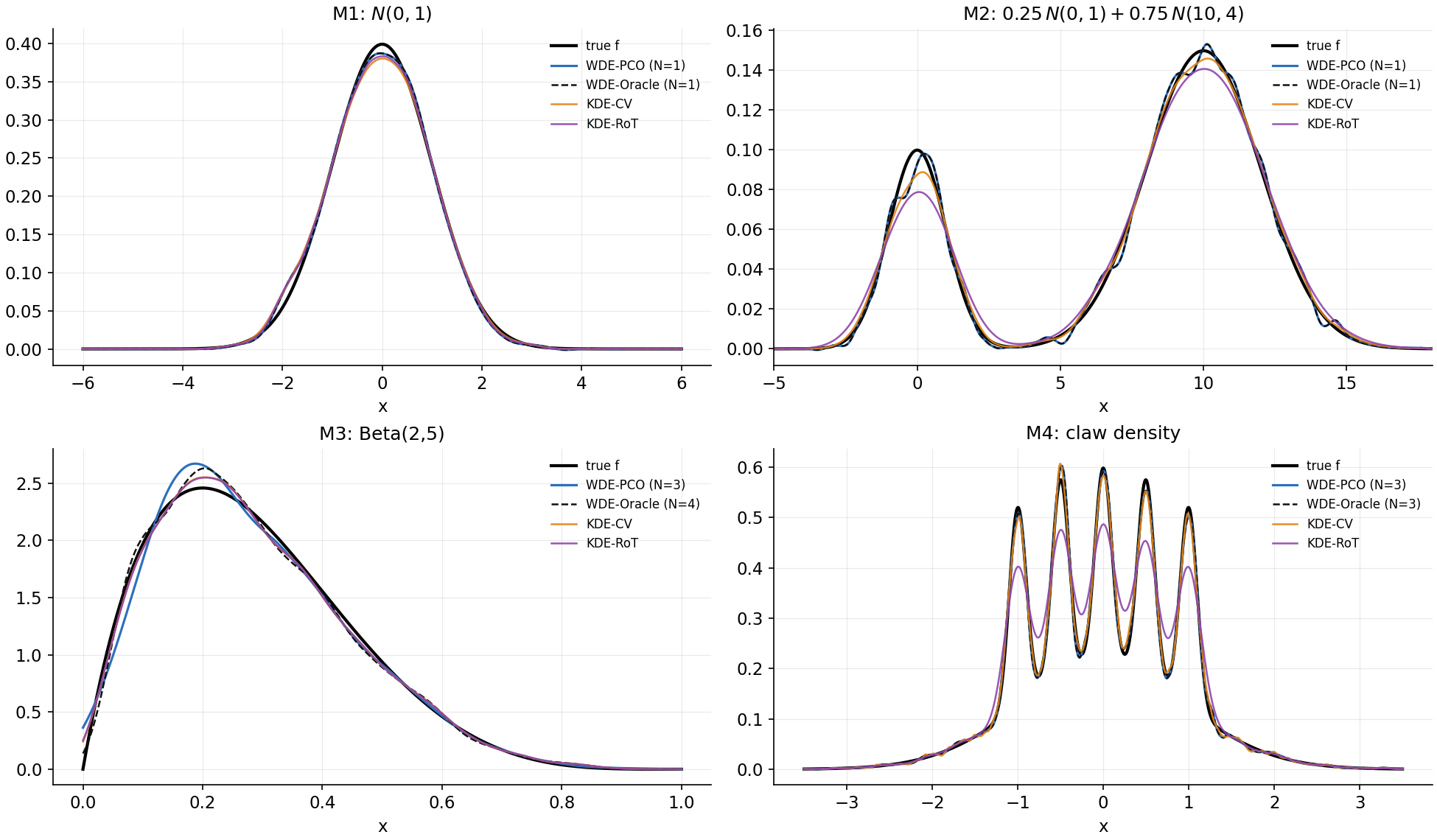}
\caption{Representative reconstructions of the {four} target
densities {: M1 and M2 (top row), M3 and M4 (bottom row)}.}
\label{fig:reconstruction}
\end{figure}

Table~\ref{tab:simulation} reports the MISE and the standard deviation of the ISE
across the $M=100$ replications for each method. Figure~\ref{fig:boxplots} shows
the corresponding boxplots of the ISE, and Figure~\ref{fig:reconstruction} shows
representative density reconstructions.

\newpage 
\begin{table}[h!]
\centering
\begin{tabular}{llrr}
\hline
Model & Method & MISE ($\times 10^{3}$) & SD ($\times 10^{3}$)\\
\hline
\multirow{4}{*}{M1} & WDE-PCO    & 0.437 & 0.277\\
                    & WDE-Oracle & 0.437 & 0.277\\
                    & KDE-CV     & 0.428 & 0.230\\
                    & KDE-RoT    & 0.402 & 0.208\\
\hline
\multirow{4}{*}{M2} & WDE-PCO    & 0.484 & 0.148\\
                    & WDE-Oracle & 0.484 & 0.148\\
                    & KDE-CV     & 0.252 & 0.088\\
                    & KDE-RoT    & 0.991 & 0.212\\
\hline
\multirow{4}{*}{M3} & WDE-PCO    & 7.211 & 1.220\\
                    & WDE-Oracle & 3.697 & 1.490\\
                    & KDE-CV     & 3.606 & 3.013\\
                    & KDE-RoT    & 2.861 & 1.092\\
\hline
{\multirow{4}{*}{M4}} & {WDE-PCO}    & {0.388}  & {0.075}\\
                    & {WDE-Oracle} & {0.387}  & {0.076}\\
                    & {KDE-CV}     & {0.458}  & {0.090}\\
                    & {KDE-RoT}    & {11.618} & {0.347}\\
\hline
\end{tabular}
{\caption{Monte Carlo estimates of the mean integrated squared error
		(MISE) and the standard deviation of the integrated squared error (ISE)
		over $M=100$ replications (in units of $10^{-3}$). M1--M3 use
		$n=2^{12}$, whereas M4 uses $n=2^{15}$. Since M4 uses a larger sample size, its absolute MISE values are not directly comparable with those of M1--M3; the comparisons across methods within each model remain meaningful.\label{tab:simulation}}}

\end{table}

The selected level reflects the structure of each target. WDE-PCO matches its
oracle on the two Gaussian targets M1 and M2:
it selects $\widehat N=1$ in every replication, exactly matching
$N_{\mathrm{Oracle}}$, so the two methods have identical risk
($\mathrm{MISE}\approx 4.4\times10^{-4}$ for M1 and $4.8\times10^{-4}$ for M2).
On the Beta density M3, WDE-PCO selects $\widehat N=3$, whereas the
oracle level is approximately $N_{\mathrm{Oracle}}=4$. As a result,
the MISE of WDE-PCO is about 1.95 times the oracle benchmark
($7.2\times10^{-3}$ versus $3.7\times10^{-3}$).
{On the claw density M4, WDE-PCO selects $\widehat N=3$ in every
replication, while the oracle level is also concentrated around $N_{\mathrm{Oracle}}=3$. The resulting MISE, $0.388\times10^{-3}$, is nearly identical to the oracle benchmark
$0.387\times10^{-3}$.}
Compared with the kernel estimators, the relative performance depends on
the target. On the smooth M1, all four methods are within approximately
10\% of one another. On the multimodal M2, WDE-PCO substantially
outperforms KDE-RoT, whose risk is about twice as large, although KDE-CV
achieves the smallest MISE. On the boundary-affected M3, both kernel
estimators outperform WDE-PCO. In contrast, on the spiky M4,
KDE-RoT strongly oversmooths the narrow peaks, while WDE-PCO achieves
the smallest MISE among the data-driven methods and remains very close
to the oracle benchmark.

\newpage 

{In conclusion, the simulations illustrate that the selected resolution
	responds to the structural complexity of the target. For the smooth
	Gaussian models M1 and M2, the PCO rule selects the coarse level
	$N=1$ and coincides with the oracle. For the more structured Beta and
	claw densities, it selects the higher level $\widehat N=3$. Its
	performance remains within a factor of two of the oracle on M3 and is
	nearly identical to the oracle benchmark on the spiky M4. At the same
	time, the weaker performance on the boundary-supported Beta density
	shows that the relative advantage of the wavelet estimator depends on
	the structure of the target. Overall, the experiments support the
	ability of the PCO rule to adapt the resolution level without selecting
	the maximal overfitting level $N_{\max}$.}

\section{Proofs}\label{sec:proofs}

The proofs rely on a control of the stochastic error of the
projection estimator at a fixed resolution level. We state it first, together
with an oracle-type upper bound and the two consequences in which the
concentration bound is used, and prove the concentration statement at
the end of the section; the proofs of Theorems~\ref{thm:oracle}
and~\ref{thm:rate} then follow. In the sequel, we use the notation of Sections~\ref{sec:construction}
and~\ref{sec:selection}, and we recall the pointwise kernel bound
\begin{equation}\label{eq:kernel-recall}
\sum_{k\in\Z}\varphi_{N,k}^2(x)=K_N(x,x)\le(2A+1)\,\Phi_0^2\,2^{N},
\qquad x\in\R,\ N\in\N,
\end{equation}
from Lemma~\ref{lem:kernel-bound}, and the variance bound
$V_N=\E\|\widehat f_N-f_N\|^2\le(2A+1)\Phi_0^2\,2^{N}/n$
from~\eqref{eq:variance-bound}. We write
\begin{equation}\label{eq:V2-def}
V_2(n,N):=(2A+1)\,\Phi_0^2\,\frac{2^{N}}{n}+\frac{\ln n}{n}
\end{equation}
for the deterministic control rate appearing below; it is nondecreasing in $N$
and of the order $2^{N}/n$ of the variance, up to the additive term
$\ln n/n$.

\begin{proposition}\label{prop:concentration}
Assume \emph{(A1)--(A2)} and $f\in\L^\infty(\R)$. {There is a universal numerical constant
$\kappa=\frac{200}{9}$} such that, for every $p\ge1$, setting
$c_1(p):=\kappa\,(1+p)\bigl(1+\|f\|_\infty\bigr)$, {one has, for all $n$ sufficiently large and every $N\in\mathcal N$,}
\begin{equation}\label{eq:concentration}
\P\!\left(\bigl\|\widehat f_N-f_N\bigr\|^2>c_1(p)\,V_2(n,N)\right)\le n^{-p}.
\end{equation}
Here $\mathcal N=\{N:2^{N}\le n/(\ln n)^2\}$ is the candidate collection
of~\eqref{eq:N-set}, so that $2^{N}\le n/(\ln n)^2$ for every $N\in\mathcal N$.
\end{proposition}

Proposition~\ref{prop:concentration} is used in the following oracle-type upper
bound, which controls the risk of the selected estimator on a large-probability
event without any restriction on the penalty term. It links the
concentration statement and Theorem~\ref{thm:oracle}, and is established in the
course of the proof of that theorem.

\begin{proposition}\label{prop:oracle-aux}
{
Assume \emph{(A1)--(A2)} and $f\in\L^\infty(\R)$. Let $p\ge1$ and  
$\theta\in(0,1)$. Let $c_{\star}(p)$ and  $\Omega_p$ be as in Lemma~\ref{lem:good-event}, while $c_{3}(p)$ and $\Omega''_p$ be as in Lemma~\ref{lem:mixedcross}. On the event $\big( \Omega_p \cap \Omega''_p \big)$, for every
$N\in\mathcal N$, 
\begin{align*}
(1-\theta)\,\bigl\|\widehat f_{\widehat N}-f\bigr\|^2
\le(1+\theta)\,\bigl\|\widehat f_N-f\bigr\|^2
&+\bigl(\pen_\lambda(N)-\pen_\lambda(\widehat N)\bigr) \\
+ & \, 2\bigl\|\widehat f_{\widehat N}-f_{\widehat N}\bigr\|^2
+
\theta c_\star(p)  
\big( V_{2}(n,\widehat{N}) + V_{2}(n,N) \big)
+
4\frac{c_{3}(p)}{\theta} \frac{\ln n}{n}  . 
\end{align*}
Consequently, for $n$ large enough, the same inequality holds with probability
at least $1-2n^{-p}$.
}\end{proposition}

The oracle analysis is based on Proposition~\ref{prop:concentration}. We apply this proposition once to obtain a high-probability event. 

Since $|\mathcal N|=N_{\max}+1\le \log_2 n+1\le 2\log_2 n$, applying a union bound introduces only a logarithmic factor. For convenience, we define
\begin{equation}\label{eq:cstar}
c_\star(p):=c_1(p+1),
\end{equation}
where $c_1(\cdot)$ is the constant appearing in Proposition~\ref{prop:concentration}. 

\begin{lem}\label{lem:good-event}
Assume \emph{(A1)--(A2)} and $f\in\L^\infty(\R)$. For every $p\ge1$, with
$c_\star(p)$ as in~\eqref{eq:cstar}, the event
\[
\Omega_p:=\Bigl\{\,\bigl\|\widehat f_M-f_M\bigr\|^2\le c_\star(p)\,V_2(n,M)
\ \text{ for all }M\in\mathcal N\,\Bigr\}
\]
satisfies $\P(\Omega_p)\ge1-n^{-p}$ for all $n$ large enough.
\end{lem}

\begin{proof}
Apply Proposition~\ref{prop:concentration} with $p+1$ in place of $p$ at each
$M\in\mathcal N$, and take the union bound over the at most $2\log_2 n$ levels:
\[
\P(\Omega_p^c)
\le\sum_{M\in\mathcal N}\P\!\left(\|\widehat f_M-f_M\|^2>c_1(p+1)V_2(n,M)\right)
\le2(\log_2 n)\,n^{-(p+1)}\le n^{-p}
\]
for $n$ large enough, since $2\log_2 n\le n$ eventually. As $N_{\max}\in
\mathcal N$, the event $\Omega_p$ controls in particular the overfitting level
$M=N_{\max}$ and the data-driven level $M=\widehat N$.
\end{proof}

The second ingredient is a cross-term identity that gives an equality instead of a Cauchy--Schwarz bound.

\begin{lem}\label{lem:crossterm}
For every $M\le N_{\max}$,
\begin{equation}\label{eq:Pres-identity}
P_M\bigl(\widehat f_{N_{\max}}-f_{N_{\max}}\bigr)=\widehat f_M-f_M,
\end{equation}
and consequently
\begin{equation}\label{eq:crossterm-identity}
\bigl\langle\widehat f_M-f_M,\ \widehat f_{N_{\max}}-f_{N_{\max}}\bigr\rangle
=\bigl\|\widehat f_M-f_M\bigr\|^2 .
\end{equation}
In particular, on the event $\Omega_p$ of Lemma~\ref{lem:good-event},
$\bigl\langle\widehat f_M-f_M,\widehat f_{N_{\max}}-f_{N_{\max}}\bigr\rangle\le
c_\star(p)\,V_2(n,M)\le c_\star(p)\,V_2(n,N_{\max})$ for every $M\in\mathcal N$.
\end{lem}

\begin{proof}
By Lemma~\ref{lem:proj-identity}, $\widehat f_M=P_M\widehat f_{N_{\max}}$ and
$f_M=P_M f_{N_{\max}}$, so~\eqref{eq:Pres-identity} follows by linearity of
$P_M$. Since $\widehat f_M-f_M\in \V_M$ and $P_M$ is the orthogonal projection
onto $\V_M$,
\begin{align*}
\bigl\langle\widehat f_M-f_M,\ \widehat f_{N_{\max}}-f_{N_{\max}}\bigr\rangle
& =\bigl\langle\widehat f_M-f_M,\ P_M(\widehat f_{N_{\max}}-f_{N_{\max}})\bigr\rangle \\
& =\bigl\langle\widehat f_M-f_M,\ \widehat f_M-f_M\bigr\rangle
=\bigl\|\widehat f_M-f_M\bigr\|^2,
\end{align*}
which is~\eqref{eq:crossterm-identity}. The last statement is then immediate
from the definition of $\Omega_p$ and the monotonicity of $V_2(n,\cdot)$.
\end{proof}

The risk analysis in Theorem~\ref{thm:rate} relies on a sharper bound for the mixed cross term. This bound is crucial for obtaining the optimal rate. The point is that this term is not the overfitting fluctuation, but a centred empirical process indexed by a deterministic function. Hence its standard deviation is of order $\|f_M-f\|/\sqrt n$, rather than $\sqrt{2^{N_{\max}}/n}$.

\begin{lem}\label{lem:mixedcross}
Assume \emph{(A1)--(A2)} and $f\in \L^\infty(\R)$. Let
\[
\mu_n(g):=\frac1n\sum_{u=1}^{n}
\bigl(g(X_u)-\E[g(X_1)]\bigr).
\]
Then, for every $M\le N_{\max}$,
\begin{equation}\label{eq:T3-identity}
\bigl\langle f_M-f,\ \widehat f_{N_{\max}}-f_{N_{\max}}\bigr\rangle
=
\mu_n\bigl(f_M-f_{N_{\max}}\bigr).
\end{equation}
Moreover, there exists a constant $C_K<\infty$, depending only on $\varphi$,
such that
\[
\|f_N\|_\infty\le C_K\|f\|_\infty,
\qquad N\in\mathcal N.
\]
Consequently, for every $p\ge1$, there exists a constant $c_3(p)$, depending only on
$p$, $\|f\|_\infty$, and $\varphi$, such that, on an event $\Omega_p''$ satisfying
\[
\P(\Omega_p'')\ge 1-n^{-p},
\]
one has, simultaneously for all $M\in\mathcal N$ and all $\eta\in(0,1)$,
\begin{equation}\label{eq:T3-bound}
2\bigl| \bigl\langle f_M-f,\ \widehat f_{N_{\max}}-f_{N_{\max}}\bigr\rangle \bigr|
\le \eta\,\|f_M-f\|^2 + \frac{c_3(p)}{\eta}\,\frac{\ln n}{n}.
\end{equation}
\end{lem}

\begin{proof}

Set $N=N_{\max}$. For any $g\in \L^2(\R)$, the reproducing property of the
projection kernel gives
\[
\langle g,K_N(x,\cdot)\rangle=(P_Ng)(x).
\]
Using the definitions of $\widehat f_N$ and $f_N=P_Nf$, we get
\[
\begin{aligned}
\bigl\langle g,\widehat f_N-f_N\bigr\rangle
&= \frac1n\sum_{u=1}^n
\bigl\langle g,K_N(X_u,\cdot)\bigr\rangle - \langle g,P_Nf\rangle  \\
&= \frac1n\sum_{u=1}^n (P_Ng)(X_u) - \langle P_Ng,f\rangle  \\
&= \mu_n(P_Ng).
\end{aligned}
\]
We now apply this identity with $g=f_M-f$. Since $M\le N_{\max}$, we have
$\V_M\subset \V_{N_{\max}}$, hence
\[
P_{N_{\max}}f_M=f_M,
\qquad
P_{N_{\max}}f=f_{N_{\max}}.
\]
Therefore
\[
P_{N_{\max}}(f_M-f)=f_M-f_{N_{\max}}.
\]
It follows that
\[
\bigl\langle f_M-f,\widehat f_{N_{\max}}-f_{N_{\max}}\bigr\rangle
=
\mu_n(f_M-f_{N_{\max}}),
\]
which proves \eqref{eq:T3-identity}.

\medskip
\noindent
Now, we establish a uniform $\L^\infty$ bound for $f_N$. By the definition of the projection kernel,
\[
K_N(x,y)=\sum_k \varphi_{N,k}(x)\varphi_{N,k}(y).
\]
Using \eqref{eq:kernel-recall} together with assumptions \emph{(A1)--(A2)}, we obtain
\[
\begin{aligned}
\int_{\R}|K_N(x,y)|\,dy
&\le
\sum_k |\varphi_{N,k}(x)|
\int_{\R}|\varphi_{N,k}(y)|\,dy  \\
&=
2^{-N/2}\|\varphi\|_1
\sum_k 2^{N/2}|\varphi(2^Nx-k)|  \\
&\le
\|\varphi\|_1(2A+1)\Phi_0
=:C_K.
\end{aligned}
\]
The constant $C_K$ is independent of $N$ and $x$. Hence
\[
\begin{aligned}
\|f_N\|_\infty
&=
\sup_x
\left|
\int_{\R}K_N(x,y)f(y)\,dy
\right|  \\
&\le
\sup_x
\int_{\R}|K_N(x,y)|\,|f(y)|\,dy  \le
C_K\|f\|_\infty.
\end{aligned}
\]

Finally, fix $M\in\mathcal N$ and write
\[
g_M:=f_M-f_{N_{\max}}.
\]
By the uniform $\L^\infty$ bound as derived above,
\[
\|g_M\|_\infty
\le
\|f_M\|_\infty+\|f_{N_{\max}}\|_\infty
\le
2C_K\|f\|_\infty.
\]
Moreover,
\[
\begin{aligned}
\Var\bigl(g_M(X_1)\bigr)
\le
\E[g_M(X_1)^2]  =
\int_{\R}g_M(x)^2 f(x)\,dx  \le
\|f\|_\infty\|g_M\|^2.
\end{aligned}
\]
Bernstein's inequality then gives the following bound: for every $x\ge1$,
\[
\P\left(
|\mu_n(g_M)|
>
C\left[
\|g_M\|\sqrt{\frac{x}{n}}
+
\frac{x}{n}
\right]
\right)
\le
2e^{-x},
\]
where $C$ depends only on $\|f\|_\infty$ and $C_K$.

Choose $x=(p+2)\ln n$. Since
\[
|\mathcal N|=N_{\max}+1\le 2\log_2 n,
\]
a union bound over $M\in\mathcal N$ gives an event $\Omega_p''$ such that
\[
\P(\Omega_p'')\ge 1-n^{-p},
\]
and on this event,
\begin{equation}\label{eq:mixed-mun-bound}
|\mu_n(g_M)|
\le
C_p\left[
\|g_M\|\sqrt{\frac{\ln n}{n}}
+
\frac{\ln n}{n}
\right],
\qquad M\in\mathcal N,
\end{equation}
where $C_p$ depends only on $p$, $\|f\|_\infty$, and $C_K$.

It remains to compare $\|g_M\|$ with the bias term $\|f_M-f\|$. Since
\[
f_M-f
=
(f_M-f_{N_{\max}})
+
(f_{N_{\max}}-f),
\]
and since $f_M-f_{N_{\max}}\in \V_{N_{\max}}$ while
$f_{N_{\max}}-f\perp \V_{N_{\max}}$, Pythagoras' identity gives
\[
\|f_M-f\|^2
=
\|f_M-f_{N_{\max}}\|^2
+
\|f_{N_{\max}}-f\|^2.
\]
Therefore
\[
\|g_M\|
=
\|f_M-f_{N_{\max}}\|
\le
\|f_M-f\|.
\]

Combining this inequality with \eqref{eq:mixed-mun-bound}, and using
\eqref{eq:T3-identity}, we obtain on $\Omega_p''$ that
\[
\begin{aligned}
2\bigl|
\bigl\langle f_M-f,\widehat f_{N_{\max}}-f_{N_{\max}}\bigr\rangle
\bigr|
&=
2|\mu_n(g_M)|  \\
&\le
2C_p\|f_M-f\|\sqrt{\frac{\ln n}{n}}
+
2C_p\frac{\ln n}{n}.
\end{aligned}
\]
Finally, for every $\eta\in(0,1)$, the inequality
$2ab\le \eta a^2+\eta^{-1}b^2$ yields
\[
2C_p\|f_M-f\|\sqrt{\frac{\ln n}{n}}
\le
\eta\|f_M-f\|^2
+
\frac{C_p^2}{\eta}\frac{\ln n}{n}.
\]
Since $\eta\in(0,1)$, we also have
\[
2C_p\frac{\ln n}{n}
\le
\frac{2C_p}{\eta}\frac{\ln n}{n}.
\]
Thus
\[
2\bigl|
\bigl\langle f_M-f,\widehat f_{N_{\max}}-f_{N_{\max}}\bigr\rangle
\bigr|
\le
\eta\|f_M-f\|^2
+
\frac{C_p^2+2C_p}{\eta}\frac{\ln n}{n}.
\]
Renaming $C_p^2+2C_p$ as $c_3(p)$ gives \eqref{eq:T3-bound}.
\end{proof}

\subsection{Proof of Proposition~\ref{prop:concentration}}

The proof rests on Talagrand's concentration inequality for the supremum of a
bounded empirical process, which we recall in the form due to Bousquet.

\begin{lem}[Talagrand--Bousquet inequality]\label{lem:talagrand}
Let $X_1,\dots,X_n$ be i.i.d.\ random variables and let $\mathcal F$ be a
countable class of real measurable functions for which there exist constants
$b,v,H>0$ such that
\[
\sup_{t\in\mathcal F}\|t\|_\infty\le b,
\qquad
\sup_{t\in\mathcal F}\Var\bigl(t(X_1)\bigr)\le v,
\qquad
\E\Bigl[\sup_{t\in\mathcal F}\nu_n(t)\Bigr]\le H,
\]
where $\nu_n(t):=\frac1n\sum_{u=1}^{n}\bigl(t(X_u)-\E[t(X_1)]\bigr)$ is the
centred empirical process indexed by $\mathcal F$.
Then, for every $x>0$ and every $\varepsilon>0$, with probability at least
$1-e^{-x}$,
\begin{equation}\label{eq:talagrand}
\sup_{t\in\mathcal F}\nu_n(t)
\le(1+\varepsilon)\,H+\sqrt{\frac{2\,v\,x}{n}}
+\Bigl(\tfrac13+\tfrac1\varepsilon\Bigr)\frac{b\,x}{n}.
\end{equation}
\end{lem}

The inequality of Lemma~\ref{lem:talagrand} is due to Bousquet
\cite{bousquet2002bennett}; we use the form stated in Massart
\cite{massart2007concentration}.

\medskip

We now prove Proposition~\ref{prop:concentration}. The goal is to control
\[
\|\widehat f_N-f_N\|^2.
\]
The key observation is that $\|\widehat f_N-f_N\|$ can be written as the
supremum of a centred empirical process over the unit ball of $\V_N$. Let $\mathcal B_N:=\{t\in \V_N:\|t\|\le 1\}$.

We claim that
\begin{equation}\label{eq:variational}
\|\widehat f_N-f_N\|
=
\sup_{t\in\mathcal B_N}\nu_n(t).
\end{equation}
Indeed, write $t=\sum_k a_k\varphi_{N,k}$, $\sum_k a_k^2\le 1$. By orthonormality,
\[
\nu_n(t)
=
\sum_k a_k(\widehat\alpha_{N,k}-\alpha_{N,k}).
\]
Taking the supremum over all sequences $(a_k)_k$ with
$\sum_k a_k^2\le1$, the equality case of the Cauchy--Schwarz inequality gives
\[
\sup_{t\in\mathcal B_N}\nu_n(t)
=
\left(
\sum_k(\widehat\alpha_{N,k}-\alpha_{N,k})^2
\right)^{1/2}
=
\|\widehat f_N-f_N\|.
\]
This proves \eqref{eq:variational}.

For measurability, the supremum in \eqref{eq:variational} may be restricted to a
countable dense subset of $\mathcal B_N$, for instance the functions with
rational coordinates and finitely many nonzero coefficients. To see this, note
that for $t\in \V_N$,
\[
|t(x)|
=
|\langle t,K_N(x,\cdot)\rangle|
\le
\|t\|\sqrt{K_N(x,x)}
\le
\|t\|\sqrt{(2A+1)\Phi_0^2\,2^N}.
\]
Thus $t\mapsto \nu_n(t)$ is continuous on $\V_N$, and the supremum over
$\mathcal B_N$ agrees with the supremum over any countable dense subset.

\medskip

We now apply Lemma~\ref{lem:talagrand} to a countable dense subset of
$\mathcal B_N$. We determine the values of $H$, $v$, and $b$.

First, by Jensen's inequality and \eqref{eq:variational},
\[
\E\left[\sup_{t\in\mathcal B_N}\nu_n(t)\right]
=
\E\|\widehat f_N-f_N\|
\le
\left(\E\|\widehat f_N-f_N\|^2\right)^{1/2}.
\]
Using the variance bound \eqref{eq:variance-bound}, we may take
\[
H
:=
\sqrt{\frac{(2A+1)\Phi_0^2\,2^N}{n}}
=
\sqrt{2A+1}\,\Phi_0\,\frac{2^{N/2}}{\sqrt n}.
\]
In particular,
\[
H^2
=
(2A+1)\Phi_0^2\, \frac{2^N}{n}.
\]

Second, for every $t\in\mathcal B_N$,
\[
\Var\bigl(t(X_1)\bigr)
\le
\E[t(X_1)^2]
=
\int_{\mathbb R} t^2(x)f(x)\,dx
\le
\|f\|_\infty\|t\|^2
\le
\|f\|_\infty.
\]
Thus we may take $v:=\|f\|_\infty$. Third, for $t=\sum_k a_k\varphi_{N,k}\in\mathcal B_N$, the Cauchy--Schwarz
inequality gives
\[
\begin{aligned}
|t(x)|
&\le
\left(\sum_k a_k^2\right)^{1/2}
\left(\sum_k \varphi_{N,k}^2(x)\right)^{1/2}  \\
&\le
\left((2A+1)\Phi_0^2\,2^N\right)^{1/2}.
\end{aligned}
\]
Hence we take
\[
b
:=
\sqrt{2A+1}\,\Phi_0\,2^{N/2}
=
H\sqrt n.
\]

\medskip
Then, we apply \eqref{eq:talagrand} with
\[
x=p\ln n,
\qquad
\varepsilon=1.
\]
Using the values of $H$, $v$, and $b$ found above, we obtain that, with
probability at least $1-n^{-p}$,
\[
\|\widehat f_N-f_N\|
\le
2H
+
\sqrt{\frac{2\|f\|_\infty p\ln n}{n}}
+
\frac43\,
\frac{\sqrt{2A+1}\,\Phi_0\,2^{N/2}p\ln n}{n}.
\]

The last term can be rewritten as
\[
\frac{4}{3}H\frac{p\ln n}{\sqrt n}.
\]
Since $N\in\mathcal N$, we have
\[
2^N\le \frac{n}{(\ln n)^2}.
\]
Therefore,
\[
H
=
\sqrt{(2A+1)\Phi_0^2\frac{2^N}{n}}
\le
\frac{\sqrt{2A+1}\Phi_0}{\ln n}.
\]
Moreover,
\[
\frac{4}{3}\sqrt{2A+1}\Phi_0 2^{N/2}
\frac{p\ln n}{n}
\le
\frac{4}{3}\sqrt{2A+1}\Phi_0
\frac{p}{\sqrt n}.
\]

For every fixed $p\ge1$, we have
\[
\frac{p\ln n}{\sqrt n}\le 1
\]
for all $n$ sufficiently large. Hence,
\[
\frac{4}{3}H\frac{p\ln n}{\sqrt n}
\le
\frac{4}{3}H.
\]
Therefore, with probability at least $1-n^{-p}$,
\[
\|\widehat f_N-f_N\|
\le
\frac{10}{3}H
+
\sqrt{2\|f\|_\infty}\,
\sqrt{\frac{p\ln n}{n}}.
\]
Thus one may take
\[
\kappa_1=\frac{10}{3},
\qquad
\kappa_2=\sqrt{2\|f\|_\infty}.
\]

\medskip

Using $(a+b)^2\le 2a^2+2b^2$, we get
\[
\begin{aligned}
\|\widehat f_N-f_N\|^2 &\le 2\kappa_1^2 H^2 + 2\kappa_2^2\,\frac{p\ln n}{n}  \\
&= 2\kappa_1^2(2A+1)\Phi_0^2\,\frac{2^N}{n} + 4\|f\|_\infty p\,\frac{\ln n}{n}.
\end{aligned}
\]
{Hence,
\[
\|\widehat f_N-f_N\|^2 \le c_1(p)V_2(n,N),
\]
where, since $\kappa_1=10/3$,
\[
2\kappa_1^2=\frac{200}{9},
\]
and therefore
\[
\max\{2\kappa_1^2,4p\|f\|_\infty\} \le \frac{200}{9}(1+p)(1+\|f\|_\infty).
\]
Thus, one may take
\[
c_1(p) = \frac{200}{9}(1+p)(1+\|f\|_\infty).
\]
This is exactly~\eqref{eq:concentration}. This completes the proof.}   \hfill$\qed$

\subsection{Proof of Proposition~\ref{prop:oracle-aux}}
{   
We have 
\[
\begin{aligned}
\bigl\|\widehat f_N-\widehat f_{N_{\max}}\bigr\|^2
&=
\bigl\|\widehat f_N-f\bigr\|^2
+
\bigl\|\widehat f_{N_{\max}}-f\bigr\|^2  
-
2\bigl\langle
\widehat f_N-f,\widehat f_{N_{\max}}-f
\bigr\rangle ,
\end{aligned}
\]
so the criterion $\Crit(N)$ can be written as
\[
\begin{aligned}
\Crit(N)
&=
\bigl\|\widehat f_N-f\bigr\|^2
-
2\bigl\langle
\widehat f_N-f,\widehat f_{N_{\max}}-f
\bigr\rangle
+
\pen_\lambda(N)  
+
\bigl\|\widehat f_{N_{\max}}-f\bigr\|^2 .
\end{aligned}
\]
Since the last term is independent of $N$,  one has that $\widehat N$ (defined in~\eqref{eq:hatN}) also minimizes
\[
N\longmapsto
\bigl\|\widehat f_N-f\bigr\|^2
-
2\bigl\langle
\widehat f_N-f,\widehat f_{N_{\max}}-f
\bigr\rangle
+
\pen_\lambda(N).
\]
Comparing the value at $\widehat N$ with the value at an arbitrary
$N\in\mathcal N$, we get
\begin{align}
\bigl\|\widehat f_{\widehat N}-f\bigr\|^2
\le\;&
\bigl\|\widehat f_N-f\bigr\|^2
+
\bigl(\pen_\lambda(N)-\pen_\lambda(\widehat N)\bigr) \notag\\
&+
2\bigl\langle
\widehat f_{\widehat N}-f,\widehat f_{N_{\max}}-f
\bigr\rangle
-
2\bigl\langle
\widehat f_N-f,\widehat f_{N_{\max}}-f
\bigr\rangle .
\label{eq:argmin}
\end{align}
Now, we decompose the cross term. For any $M\in\mathcal N$ fixed, we write
\[
\widehat f_M-f
=
(\widehat f_M-f_M)+(f_M-f)
\]
and
\[
\widehat f_{N_{\max}}-f = (\widehat f_{N_{\max}}-f_{N_{\max}}) + (f_{N_{\max}}-f).
\]
This gives
\begin{align}
\bigl\langle
\widehat f_M-f,\widehat f_{N_{\max}}-f
\bigr\rangle
&=
\bigl\langle
\widehat f_M-f_M,\widehat f_{N_{\max}}-f_{N_{\max}}
\bigr\rangle \notag
+
\bigl\langle
\widehat f_M-f_M,f_{N_{\max}}-f
\bigr\rangle \notag\\
&\quad
+
\bigl\langle
f_M-f,\widehat f_{N_{\max}}-f_{N_{\max}}
\bigr\rangle \notag
+
\bigl\langle
f_M-f,f_{N_{\max}}-f
\bigr\rangle \notag\\
&=:T_1(M)+T_2(M)+T_3(M)+T_4(M).
\label{eq:fourterms}
\end{align}
First, by the orthogonality
\begin{equation}\label{eq:ortho}
f_{N_{\max}}-f = -(I-P_{N_{\max}})f \perp \V_{N_{\max}}
\supseteq \V_M,
\qquad M\le N_{\max}, 
\end{equation}
this implies that $T_2(M) = \bigl\langle
\widehat f_M-f_M,f_{N_{\max}}-f
\bigr\rangle = 0$ since  
\[
\big( \widehat f_M-f_M \big) \in \V_M\subseteq \V_{N_{\max}},
\textrm{ while } \big( f_{N_{\max}}-f \big)  \perp \V_{N_{\max}}. 
\]
Second, $T_4(M)$ is independent of $M$. Indeed, since 
\[
f_M-f
=
(f_M-f_{N_{\max}})+(f_{N_{\max}}-f),
\]
with $\big( f_M-f_{N_{\max}} \big) \in \V_{N_{\max}}$ and
$\big( f_{N_{\max}}-f \big) \perp \V_{N_{\max}}$, we obtain  
\[
\begin{aligned}
T_4(M)
&=
\bigl\langle f_M-f,f_{N_{\max}}-f\bigr\rangle =
\bigl\langle f_{N_{\max}}-f,f_{N_{\max}}-f\bigr\rangle =
\bigl\|f_{N_{\max}}-f\bigr\|^2.
\end{aligned}
\]
Thus, $T_4(\widehat N)=T_4(N)$, so the $T_4$ terms cancel in \eqref{eq:argmin}.
\leavevmode \\  
Hence, substituting \eqref{eq:fourterms} into \eqref{eq:argmin}, and using
$T_2\equiv0$ together with the cancellation of $T_4$, we get
\begin{align}
\bigl\|\widehat f_{\widehat N}-f\bigr\|^2
\le\;&
\bigl\|\widehat f_N-f\bigr\|^2
+
\bigl(\pen_\lambda(N)-\pen_\lambda(\widehat N)\bigr)  \nonumber  
\\
&+
2T_1(\widehat N)-2T_1(N)
+
2T_3(\widehat N)-2T_3(N).  
\label{eq:rate-basic1}
\end{align}
\leavevmode \\ 
Now, we control the $T_{3}(\cdot)$ terms as follows. By Lemma~\ref{lem:mixedcross}, on $\Omega_p''$, for every $M\in\mathcal N$ and $\eta \in (0,1)$, 
\begin{equation}\label{eq:mixed-rate1}
	2|T_3(M)|
	\le
	\eta\|f_M-f\|^2
	+
	\frac{c_3(p)}{\eta}\frac{\ln n}{n}, 
\end{equation}
with $c_{3}(p)$ defined in Lemma~\ref{lem:mixedcross}. 
Using the triangle inequality and Lemma~\ref{lem:good-event} to get on $\Omega_p$ that  
\begin{align}
	\|f_M-f\|^2
	&\le
	2\bigl\|\widehat f_M-f\bigr\|^2
	+
	2\bigl\|\widehat f_M-f_M\bigr\|^2 
	\nonumber \\
	&\le
	2\bigl\|\widehat f_M-f\bigr\|^2
	+
	2c_\star(p)V_2(n,M).
	\label{eq:mixed-rate.2}
\end{align}
Thus, we get from \eqref{eq:mixed-rate1} and~\eqref{eq:mixed-rate.2} that on $\big( \Omega_p \cap \Omega_p'' \big)$, for any $\eta \in (0,1)$, 
\begin{align}
	2 T_3(\widehat N)
	&\le
	2\eta\bigl\|\widehat f_{\widehat N}-f\bigr\|^2
	+
	2\eta c_\star(p)V_2(n,\widehat N)
	+
	\frac{c_3(p)}{\eta}\frac{\ln n}{n},
	\label{eq:T3-hat-rate1}
	\\
	2|T_3(N)|
	&\le
	2\eta\bigl\|\widehat f_N-f\bigr\|^2
	+
	2\eta c_\star(p)V_2(n,N)
	+
	\frac{c_3(p)}{\eta}\frac{\ln n}{n}.
	\label{eq:T3-N-rate1}
\end{align}
\noindent  Moreover, by~\eqref{eq:crossterm-identity} in Lemma~\ref{lem:crossterm}, on the event $\Omega_p$, for every $ M \in\mathcal N$, 
\begin{align}   
	T_1(M)
	=
	\bigl\langle
	\widehat f_{M}-f_{M},
	\widehat f_{N_{\max}}-f_{N_{\max}}
	\bigr\rangle
	=
	\bigl\|\widehat f_{M}-f_{M}\bigr\|^2  
	\label{eq:T1.equivalence}
\end{align}   
Eventually, combining  
\eqref{eq:T3-hat-rate1}, \eqref{eq:T3-N-rate1}, \eqref{eq:T1.equivalence}
(here we keep the term $2T_1(\widehat N)=2\|\widehat
f_{\widehat N}-f_{\widehat N}\|^2$, to be absorbed later by the penalty at level
$\widehat N$, and simply discard the nonpositive term
$-2T_1(N)=-2\|\widehat f_N-f_N\|^2\le 0$), we
obtain from~\eqref{eq:rate-basic1}, on $\big( \Omega_p \cap \Omega_p'' \big)$, for every $N\in\mathcal N$ and $\eta \in (0,1)$,
\[
\begin{aligned}
	\bigl\|\widehat f_{\widehat N}-f\bigr\|^2
	\le\;&
	\bigl\|\widehat f_N-f\bigr\|^2
	+
	\bigl(\pen_\lambda(N)-\pen_\lambda(\widehat N)\bigr)   
	+
	2\eta\bigl\|\widehat f_{\widehat N}-f\bigr\|^2 
	+
	2\eta\bigl\|\widehat f_N-f\bigr\|^2 \\
	&+
	2\bigl\|\widehat f_{\widehat{N}}-f_{\widehat{N}}\bigr\|^2
	+ 
	2\eta c_\star(p)  
	\big( V_{2}(n,\widehat{N}) + V_{2}(n,N) \big)
	+
	2\frac{c_{3}(p)}{\eta} \frac{\ln n}{n}.
\end{aligned}
\]
Hence, 
\begin{align}
	(1-2\eta)
	\bigl\|\widehat f_{\widehat N}-f\bigr\|^2
	&\le
	(1+2\eta)
	\bigl\|\widehat f_N-f\bigr\|^2
	+
	\bigl(\pen_\lambda(N)-\pen_\lambda(\widehat N)\bigr)  
	\nonumber  \\ 
	& \qquad +
	2\bigl\|\widehat f_{\widehat{N}}-f_{\widehat{N}}\bigr\|^2 +
	2\eta  c_\star(p)  
	\big( V_{2}(n,\widehat{N}) + V_{2}(n,N) \big)
	+
	2\frac{c_{3}(p)}{\eta} \frac{\ln n}{n}  .  
	\label{eq:rate-good-event}
\end{align}
Choose $\eta=\theta/2$, so  $1-2\eta=1-\theta$ and
$1+2\eta=1+\theta$. 
Thus, for every $N\in\mathcal N$,
\begin{align}    
(1-\theta)\bigl\|\widehat f_{\widehat N}-f\bigr\|^2
&\le
(1+\theta)\bigl\|\widehat f_N-f\bigr\|^2
+
\bigl(\pen_\lambda(N)-\pen_\lambda(\widehat N)\bigr)  
\nonumber  \\ 
& \qquad +
2\bigl\|\widehat f_{\widehat{N}}-f_{\widehat{N}}\bigr\|^2  +
\theta c_\star(p)  
\big( V_{2}(n,\widehat{N}) + V_{2}(n,N) \big)
+
4\frac{c_{3}(p)}{\theta} \frac{\ln n}{n}  .  
\label{eq:prop62-final}
\end{align}
\noindent  
Finally, since all estimates hold on $(\Omega_p \cap \Omega_p'')$  and 
$\P(\Omega_p \cap \Omega_p'')\ge 1-2n^{-p}$ for $n$ large enough, the proof is complete. \hfill\qed
}

\subsection{Proof of Theorem~\ref{thm:oracle}}
{   
From Proposition~\ref{prop:oracle-aux}, there exists an event $\Omega: = \big( \Omega_p \cap \Omega_p'' \big)$ with
probability at least $1-2n^{-p}$, for $n$ large enough, such that on $\Omega$,
for every $N\in\mathcal N$,
\begin{align}
(1-\theta)\bigl\|\widehat f_{\widehat N}-f\bigr\|^2
&\le
(1+\theta)\bigl\|\widehat f_N-f\bigr\|^2
+
\bigl(\pen_\lambda(N)-\pen_\lambda(\widehat N)\bigr)
\nonumber  \\ 
& \qquad + 2 \bigl\|\widehat f_{\widehat N}-f_{\widehat N}\bigr\|^2
+
\theta c_\star(p)  
\big( V_{2}(n,\widehat{N}) + V_{2}(n,N) \big)
+
4\frac{c_{3}(p)}{\theta} \frac{\ln n}{n}  .
\label{eq:from-prop62}
\end{align}
It remains to control the penalty difference. Let $C_\varphi := (2A+1) \Phi_{0}^{2}$. By definition, 
\[
\pen_\lambda(M)
=
\frac{C_\varphi}{n}
\left[(\lambda+1)2^M-2^{N_{\max}}\right]
\]
hence, 
\begin{equation}\label{eq:pen-diff-rate1}
	\pen_\lambda(N)-\pen_\lambda(\widehat N)
	=
	(\lambda+1)C_\varphi\frac{2^N}{n}
	-
	(\lambda+1)C_\varphi\frac{2^{\widehat N}}{n}.
\end{equation}
On the other hand, the two $V_2(n,\cdot)$ terms in~\eqref{eq:from-prop62} satisfy
\begin{align}
	\theta c_\star(p)V_2(n,\widehat N)
	+
	\theta c_\star(p)V_2(n,N)
	\le
	2\theta c_\star(p)C_\varphi
	\frac{2^{\widehat N}+2^N}{n}
	+
	2\theta c_\star(p)\frac{\ln n}{n}.
	\label{prof:thm4.1:V2.control}
\end{align}
\indent  
Combining~\eqref{eq:pen-diff-rate1}  and~\eqref{prof:thm4.1:V2.control}, we get from~\eqref{eq:from-prop62} that on $\Omega$,  
\begin{align}
	(1-\theta)\bigl\|\widehat f_{\widehat N}-f\bigr\|^2
	&\le
	(1+\theta)\bigl\|\widehat f_N-f\bigr\|^2
	+
	\big[ (\lambda+1) + 2\theta c_{\star}(p) \big] C_\varphi\frac{2^N}{n}  
	\nonumber  \\ 
	& \qquad + 2 \bigl\|\widehat f_{\widehat N}-f_{\widehat N}\bigr\|^2 - \big[ (\lambda+1) - 2\theta c_{\star}(p) \big] C_\varphi\frac{2^{\widehat{N}}}{n}
	+
	4\frac{c_{3}(p)}{\theta} \frac{\ln n}{n}  .
	\label{eq:from-prop62:eq2}
\end{align} 
Moreover, using Lemma~\ref{lem:good-event} and the definition of $V_2(n,\cdot)$ to get that on $\Omega \subseteq \Omega_p$,
\begin{align*}
	2 \bigl\|\widehat f_{\widehat N}-f_{\widehat N}\bigr\|^2
	\le
	2 c_\star(p)V_2(n,\widehat N)  
	\le
	4c_\star(p)C_\varphi\frac{2^{\widehat N}}{n}
	+
	2c_\star(p)\frac{\ln n}{n}.
\end{align*}
Consequently, if $\lambda+1\ge (4+2\theta ) c_\star(p) $, then 
\begin{align*}
	&2\bigl\|\widehat f_{\widehat N}-f_{\widehat N}\bigr\|^2
	-
	\big[ (\lambda+1) - 2\theta c_{\star}(p) \big] C_\varphi\frac{2^{\widehat N}}{n} \notag\\
	&\qquad\le
	\bigl( (4+2\theta ) c_\star(p)-(\lambda+1)\bigr)
	C_\varphi\frac{2^{\widehat N}}{n}
	+
	2c_\star(p)\frac{\ln n}{n} \notag\\
	&\qquad\le
	2c_\star(p)\frac{\ln n}{n}.
\end{align*}
Therefore, we
obtain, on $\Omega$, for every $N\in\mathcal N$,
\begin{equation}\label{eq:rate-good-event1}
	(1-\theta) \bigl\|\widehat f_{\widehat N}-f\bigr\|^2
	\le\;
	(1+\theta) \bigl\|\widehat f_N-f\bigr\|^2
	+
	\bigl[(\lambda+1)+ 2\theta c_\star(p)\bigr]
	C_\varphi\frac{2^N}{n}
	+
	\Big( 2 c_{\star}(p) + 4 \frac{c_{3}(p)}{\theta} \Big) \frac{\ln n}{n}  .
\end{equation}
Finally, dividing both sides by $1-\theta$, we obtain
\begin{align*}
\bigl\|\widehat f_{\widehat N}-f\bigr\|^2
\le\;&
\frac{1+\theta}{1-\theta}
\bigl\|\widehat f_N-f\bigr\|^2  + 
\frac{1}{1-\theta}
\left[
\bigl[(\lambda+1)+ 2\theta c_\star(p)\bigr]
C_\varphi\frac{2^N}{n}
+
\Big( 2 c_{\star}(p) + 4 \frac{c_{3}(p)}{\theta} \Big) \frac{\ln n}{n}
\right]  
\\
\le\;&
\frac{C_{\theta,p, \lambda}}{1-\theta} \Big(
\bigl\|\widehat f_N-f\bigr\|^2 + \frac{2^N}{n} 
\Big) +
\frac{1}{1 - \theta}\Big( 2 c_{\star}(p) + 4 \frac{c_{3}(p)}{\theta} \Big) \frac{\ln n}{n}  , 
\end{align*}
where $C_{\theta,p, \lambda} := \max \big\{ 1 + \theta ;  \bigl[(\lambda+1)+ 2\theta c_\star(p)\bigr]
C_\varphi \big\}$. 
\leavevmode \\  
Since this inequality holds for every $N\in\mathcal N$, taking the minimum over
$N\in\mathcal N$ on the right-hand side yields the oracle inequality stated in
Theorem~\ref{thm:oracle}. This completes the proof.
\hfill$\qed$
}

\subsection{Proof of Theorem~\ref{thm:rate}}
 
Throughout the proof, write $C_\varphi:=(2A+1)\Phi_0^2$. 
Fix $U>0$ and let $f$ be a probability density with $\|f\|_\infty\le U$. The density $f$ enters
Proposition~\ref{prop:concentration}, Lemma~\ref{lem:good-event}, and
Lemma~\ref{lem:mixedcross} only through upper bounds involving
$\|f\|_\infty$: the variance proxy $v=\|f\|_\infty$ in the
Talagrand--Bousquet inequality, and the bounds $\|g_M\|_\infty\le 2C_K\|f\|_\infty$, $\Var(g_M(X_1))
\le \|f\|_\infty\|g_M\|^2$
in Bernstein's inequality. These bounds remain valid with
$\|f\|_\infty$ replaced by $U$. Accordingly, set $c_1(p,U):=\frac{200}{9}(p+1)(1+U)$
so that $c_1(p+1,U)=c_\star(p,U)$, with $c_\star(p,U)$ as defined before
Theorem~\ref{thm:rate}. Then Proposition~\ref{prop:concentration} and
Lemma~\ref{lem:good-event} hold uniformly over all densities satisfying
$\|f\|_\infty\le U$, with $c_1(p,U)$ and $c_\star(p,U)$ in place of
$c_1(p)$ and $c_\star(p)$, respectively. Likewise,
Lemma~\ref{lem:mixedcross} holds with a constant $c_3(p,U)$ depending
only on $p$, $U$, and $C_K$.

The preceding bounds hold once $n$ is large enough so that $(p+1)\ln n\le\sqrt n$ and $2\log_2 n\le n$. These conditions depend only on $p$. Hence we may choose
$n_0\ge3$, depending only on $p$, such that they are satisfied for all
$n\ge n_0$. For the remainder of this proof, we fix $n\ge n_0$ and,
to keep the notation unchanged in the estimates below, write
$c_\star(p)$ and $c_3(p)$ for the uniform constants
$c_\star(p,U)$ and $c_3(p,U)$, respectively. We also let
$\Omega_p$ and $\Omega_p''$ denote the corresponding events.

We work on
the event
\[
\Omega:=\Omega_p\cap\Omega_p'',
\]
where $\Omega_p$ is the event of Lemma~\ref{lem:good-event} and
$\Omega_p''$ is the event of Lemma~\ref{lem:mixedcross}. \\
By a union bound,
\[
\P(\Omega)\ge 1-2n^{-p}.
\]
All estimates below are deterministic on $\Omega$. 
\leavevmode \\  
As in the proof of Theorem~\ref{thm:oracle}, from~\eqref{eq:rate-good-event1}, we get
\begin{equation}\label{eq:rate-good-event2}
(1-\theta) \bigl\|\widehat f_{\widehat N}-f\bigr\|^2
\le\;
(1+\theta) \bigl\|\widehat f_N-f\bigr\|^2
+
\bigl[(\lambda+1)+ 2\theta c_\star(p)\bigr]
C_\varphi\frac{2^N}{n}
+
\Big( 2 c_{\star}(p) + 4 \frac{c_{3}(p)}{\theta} \Big) \frac{\ln n}{n}  .
\end{equation}
We now choose $\theta=1/2$. Then $1-\theta=1/2$, and, with the
above convention, the condition on the penalty becomes $\lambda\ge 5c_\star(p)-1=\lambda_0(p,U)$.

Under this condition, the term involving $c_\star(p)$ in the leading variance
coefficient is bounded by a constant depending only on $\lambda$. Thus the
leading constant may be chosen to depend only on $\lambda$.

For any $N\in\mathcal N$ fixed, taking the  expectation in~\eqref{eq:rate-good-event2} over the event $\Omega$, and using
\[
\E\bigl\|\widehat f_N-f\bigr\|^2
=
\|f_N-f\|^2+V_N
\le
\|f_N-f\|^2
+
C_\varphi\frac{2^N}{n},
\]
we obtain
\begin{equation}\label{eq:rate-oracle-exp}
\E\left[
\bigl\|\widehat f_{\widehat N}-f\bigr\|^2
\mathbf 1_{\Omega}
\right]
\le
C_1
\left(
\|f_N-f\|^2
+
C_\varphi\frac{2^N}{n}
\right)
+
C_2\frac{\ln n}{n},
\qquad N\in\mathcal N,
\end{equation}
where $C_1$ depends only on $\lambda$, while $C_2$ depends only on {\color{blue}$p$, $U$ and $\varphi$}.

\medskip
It remains to control
\[
\E\left[
\bigl\|\widehat f_{\widehat N}-f\bigr\|^2
\mathbf 1_{\Omega^c}
\right].
\]
For every $N\in\mathcal N$, by Cauchy--Schwarz,
\[
\widehat\alpha_{N,k}^2
=
\left(
\frac1n\sum_{u=1}^n \varphi_{N,k}(X_u)
\right)^2
\le
\frac1n\sum_{u=1}^n \varphi_{N,k}^2(X_u).
\]
Therefore, using \eqref{eq:kernel-recall},
\[
\begin{aligned}
\|\widehat f_N\|^2
=
\sum_k \widehat\alpha_{N,k}^2 \le
\frac1n\sum_{u=1}^n \sum_k \varphi_{N,k}^2(X_u) \le
C_\varphi 2^N
\le
C_\varphi n.
\end{aligned}
\]
Since $\widehat N\in\mathcal N$, the same bound holds for
$\widehat f_{\widehat N}$. Hence
\[
\bigl\|\widehat f_{\widehat N}-f\bigr\|^2
\le
2\|\widehat f_{\widehat N}\|^2+2\|f\|^2
\le
C_4 n
\]
for a constant $C_4$ depending only on $C_\varphi$ and $U$, since $\|f\|^2\le\|f\|_\infty\int_\R f\le U$.

Since $\P(\Omega^c)\le 2n^{-p}$, for $p\ge2$,
\[
\E\!\left[
\|\widehat f_{\widehat N}-f\|^2\mathbf 1_{\Omega^c}
\right]
\le
2C_4 n^{1-p}
\le
\frac{2C_4}{n}
\le
2C_4\frac{\ln n}{n}.
\]

Combining this estimate with \eqref{eq:rate-oracle-exp} and replacing $C_2$ by $C_2+2C_4$, which still depends only on $p$, $U$ and $\varphi$, we obtain
\[
\E\bigl\|\widehat f_{\widehat N}-f\bigr\|^2
\le
C_1
\left(
\|f_N-f\|^2
+
C_\varphi\frac{2^N}{n}
\right)
+
C_2\frac{\ln n}{n},
\qquad N\in\mathcal N.
\]
Taking the minimum over $N\in\mathcal N$ proves the oracle inequality
\eqref{eq:risk-oracle}.

\medskip
For the convergence rate, fix $0<r<\mathcal{S}$, $1\le q\le\infty$ and $L>0$,
and let $f\in\mathcal D_{r,q}(L,U)$. Set
\[
N_*:=\Bigl\lfloor\frac{\log_2 n}{2r+1}\Bigr\rfloor,
\qquad\text{so that}\qquad
\frac12\,n^{1/(2r+1)}<2^{N_*}\le n^{1/(2r+1)} .
\]
There is an integer $n_1\ge n_0$, depending only on $p$ and $r$, such that for
$n\ge n_1$ we have $N_*\ge1$, $n^{1/(2r+1)}\le n/(\ln n)^2$ (hence
$N_*\in\mathcal N$) and $\ln n\le n^{1/(2r+1)}$. By
Proposition~\ref{prop:bias} and $\|f\|_{B^r_{2,q}}\le L$,
\[
\|f_{N_*}-f\|^2\le C_{\mathrm{Besov}}L^2\,2^{-2rN_*}
\le4^rC_{\mathrm{Besov}}L^2\,n^{-\frac{2r}{2r+1}},
\qquad
C_\varphi\frac{2^{N_*}}{n}\le C_\varphi\,n^{-\frac{2r}{2r+1}},
\]
and $\ln n/n\le n^{-2r/(2r+1)}$. 

Applying~\eqref{eq:risk-oracle} with
$N=N_*$ gives, for $n\ge n_1$,
\[
\E\bigl\|\widehat f_{\widehat N}-f\bigr\|^2
\le\Bigl(C_1\bigl(4^rC_{\mathrm{Besov}}L^2+C_\varphi\bigr)+C_2\Bigr)
n^{-\frac{2r}{2r+1}} .
\]

The constant on the right-hand side depends only on $\lambda$, $p$, $U$,
$L$, $r$ and $\varphi$, and not on $f$. Taking the supremum over
$f\in\mathcal D_{r,q}(L,U)$ proves~\eqref{eq:rate}.  \hfill\qed

\subsubsection*{Acknowledgement}
We would like to express our sincere gratitude to the anonymous referee for the careful reading and thoughtful, constructive comments, which have significantly improved the quality, rigor, and presentation of this paper. This research is funded by Vietnam National University, Ho Chi Minh City (VNU-HCM), under grant number C2022-18-02.

\end{document}